\documentclass[hidelinks]{article}
\usepackage{amsmath}
\usepackage{color}
\usepackage{amsfonts}
\usepackage{amsmath}
\usepackage{amssymb}
\usepackage{hyperref}
\usepackage{ulem}
\usepackage{pstricks}  
\usepackage{graphicx}

\newtheorem{Thm}{Theorem}

\newtheorem{Cor}{Corollary}  
\newtheorem{lemma}{Lemma}

\newenvironment{dem}{\noindent \textit{Proof:} }{\quad \hfill $\square$}
\usepackage[left=2cm,top=2cm,right=2cm,bottom=2cm]{geometry}
\usepackage{graphicx}
\usepackage{enumitem}
\newcommand{\etalchar}[1]{$^{#1}$}

\newcommand{\Sy}{{\mathfrak S}}

\newcommand{\PP}{{\mathcal P}}
\newcommand{\G}{{\mathfrak G}}

\newcommand{\T}{{\mathfrak t}}
\newcommand{\s}{{\mathfrak s}}
\newcommand{\U}{{\mathfrak u}}
\newcommand{\V}{{\mathfrak v}}
\newcommand{\g}{{\mathfrak g}}

\newcommand{\Par}{\operatorname{Par}}
\newcommand{\Comp}{\operatorname{Comp}}
\newcommand{\rad}{\operatorname{rad}}
\newcommand{\Cont}{{\cal C}}
\newcommand{\Sym}{\operatorname{Sym}}

\newcommand{\Span}{{\rm span}}

\newcommand{\Res}{{\rm res}}
\newcommand{\std}{{\rm Std}}
\newcommand{\tab}{{\rm Tab}}

\newcommand{\C}{{\mathbb C}}
\newcommand{\Z}{{\mathbb Z}}
\newcommand{\Q}{{\mathbb Q}}

\newcommand{\F}{{\mathbb  F}}

\edef\savecatcodeat{\the\catcode`@}
\catcode`\@=11

\def\tb@ifSpecChars#1#2{#1}
\def\tb@ifNoSpecChars#1#2{#2}

\def\tableau{%
  \bgroup% matched in \tb@tableauD
  \@ifstar{\let\Tif\tb@ifNoSpecChars\tb@tableauB}% *, don't use special chars
          {\let\Tif\tb@ifSpecChars\tb@tableauB}}% no *, use special chars

\def\tb@tableauB{% add [] if no [options]
  \@ifnextchar[{\tb@tableauC}{\tb@tableauC[]}}

\def\tb@tableauC[#1]{\hbox\bgroup%
    \let\\=\cr% end line
    \def\bl{\global\let\tbcellF\tb@cellNF}%
    \def\tf{\global\let\tbcellF\tb@cellH}% highlighted cell
    \dimen2=\ht\strutbox \advance\dimen2 by\dp\strutbox%
    \ifx\baselinestretch\undefined\relax%
    \else%
       \dimen0=100sp \dimen0=\baselinestretch\dimen0%
       \dimen2=100\dimen2 \divide\dimen2 by\dimen0%
    \fi%
    \let\tpos\tb@vcenter% default position
    \tb@initYoung% default tableau type
    \tb@options#1\eoo% parse options
    \let\arrow\tb@arrow%
    \dimen0=\Tscale\dimen2%
    \dimen1=\dimen0 \advance\dimen1 by \tb@fframe%
    \lineskip=0pt\baselineskip=0pt% line spacing will be from \vbox to \dimen0
    \def\tb@nothing{}%
    \def\endcellno{$\rss\egroup\bss\egroup}% end cell w/o overlap
    \def\endcell{\endcellno\kern-\dimen0}% end cell & prepare to overlap it
    \def\begincell{\vbox to\dimen0\bgroup\vss\hbox to\dimen0\bgroup\hss$}%
    \let\overlay\tb@overlay%
    \let\fl\tb@fl%
    \let\lss\hss\let\rss\hss\let\tss\vss\let\bss\vss% cell alignment
    \def\mkcell##1{% format individual cell
        \let\tbcellF\tb@cellD% default cell frame
        \def\tb@cellarg{##1}% store cell contents
        \ifx\tb@cellarg\tb@nothing\let\tb@cellarg\tb@cellE\fi%
        \begincell\tb@cellarg\endcellno% the actual cell content
        \tbcellF}% draw cell frame
    \let\savecellF\tbcellF% save global value of cellF in case of nested tableau
     \Tif{\catcode`,=4\catcode`|=\active}{}\tb@tableauD}%

\let\tb@savetableauD\tableauD% save any current definition
{% set up characters which will be interpreted as command characters
    \catcode`|=\active \catcode`*=\active \catcode`~=\active%
    \catcode`@=\active% command characters
\gdef\tableauD#1{%
  \Tif{% make all the command characters active in math mode when #1 parsed
    \mathcode`|="8000 \mathcode`*="8000%
    \mathcode`~="8000 \mathcode`@="8000%
    \def@{\bullet}%
    \let|\cr% end line
    \let*\tf% highlighted cell
    \let~\sk% skew cell
  }{}%
  \tpos{\tabskip=0pt\halign{&\mkcell{##}\cr#1\crcr}}%
  \global\let\tbcellF\savecellF% restore global value
  \egroup% match \hbox\bgroup at start of \tableauC
  \egroup}% match \bgroup at start of \tableau
}
\let\tb@tableauD\tableauD% rename the command
\let\tableauD\tb@savetableauD% restore old command with this name
\let\tb@savetableauD\undefined

\def\tb@options#1{\ifx#1\eoo\relax\else\tb@option#1\expandafter\tb@options\fi}

\def\tb@option#1{%
  \if#1t\let\tpos\tb@vtop\fi%        t = align at top
  \if#1c\let\tpos\tb@vcenter\fi%     c = align at center
  \if#1b\let\tpos\vbox\fi%           b = align at bottom
  \if#1F\tb@initFerrers\fi%          F = Ferrers diagram
  \if#1Y\tb@initYoung\fi%            Y = Young diagram
  \if#1s\tb@initSmall\fi%            s = small boxes
  \if#1m\tb@initMedium\fi%           m = medium boxes
  \if#1l\tb@initLarge\fi%            l = large boxes
  \if#1p\tb@initPartition\fi%            p = small partition sized boxes
  \if#1a\tb@initArrow\fi%            a = use arrow font as base dimension
}

\def\tb@vcenter#1{\ifmmode\vcenter{#1}\else$\vcenter{#1}$\fi}

\def\tb@vtop#1{\hbox{\raise\ht\strutbox\hbox{\lower\dimen0\vtop{#1}}}}

\def\tb@initPartition{\def\Tscale{.3}}
\def\tb@initSmall{\def\Tscale{1}}
\def\tb@initMedium{\def\Tscale{2}}
\def\tb@initLarge{\def\Tscale{3}}

\def\tb@initArrow{\dimen2=1.25em}

\def\tb@initYoung{%
  \def\tb@cellE{}% empty cell stays empty
  \let\tb@cellD\tb@cellN% default frame is normal frame
  \def\sk{\global\let\tbcellF\tb@cellNF}}% skew cells are empty
\def\tb@initFerrers{%
  \def\tb@cellE{\bullet}% empty cell gets bullet
  \let\tb@cellD\tb@cellNF% default frame is no frame
  \def\sk{\bullet}}% skew cell gets bullet

\tb@initMedium% default scale

\def\tb@sframe#1{%
  \vbox to0pt{%            Embed frame in a box of no vert or hor extent
    \vss%                            pull box above reference point
    \hbox to0pt{%
      \hss%                          pull box left of reference point
      \vbox to\dimen1{%              Actual width of frame
        \hrule depth #1 height0pt% draw top edge of frame
        \vss%                     begin vcenter sides
        \hbox to\dimen1{%           horiz box with side edges just inside
          \vrule width #1 height\dimen1% left edge
          \hss%                     stretch center
          \vrule width #1%         right edge
          }%
        \vss%                     end vcenter sides
        \hrule height #1 depth 0in% bottom edge
        }%
      \kern-\tb@hframe%           horiz alignment off by half line width
      }%
    \kern-\tb@hframe}}%           vert alignment off by half line width
\def\tb@hframe{.2pt}\def\tb@fframe{.4pt}\def\tb@bframe{1.2pt}
\def\tb@cellH{\tb@sframe{\tb@bframe}}       % bold frame
\def\tb@cellNF{}                            % no frame
\def\tb@cellN{\tb@sframe{\tb@fframe}}       % normal frame
\let\tbcellF\tb@cellN                       % default is normal

\def\tb@rpad{1pt}
\def\tb@lpad{1pt}
\def\tb@tpad{1.8pt}
\def\tb@bpad{1.8pt}

\def\tb@overlay{\endcell\@ifnextchar[{\tb@overlaya}{\begincell}}
\def\tb@overlaya[#1]{\vbox to\dimen0\bgroup%
  \tb@overlayoptions#1\eoo%
  \tss\hbox to\dimen0\bgroup\lss}
\def\tb@overlayoptions#1{\ifx#1\eoo\relax\else\tb@overlayoption#1\expandafter\tb@overlayoptions\fi}

\def\tb@overlayoption#1{
  \if#1t\def\tss{\vskip\tb@tpad}\let\bss\vss\fi% t = align at top
  \if#1c\let\tss\vss\let\bss\vss\fi%             c = align at center
  \if#1b\def\bss{\vskip\tb@bpad}\let\tss\vss\fi% b = align at bottom
  \if#1l\def\lss{\hskip\tb@lpad}\let\rss\hss\fi% l = align at left
  \if#1m\let\lss\hss\let\rss\hss\fi%             m = align at middle
  \if#1r\def\rss{\hskip\tb@rpad}\let\lss\hss\fi% r = align at right
}

\def\tb@fl{\endcell\begincell\vrule depth 0pt width \dimen0 height \dimen0 \endcell\begincell}

\@ifundefined{diagram}{}{
\def\tb@arrowpad{.5}

\newoptcommand{\tb@arrow}{\@ne}[2]{%
  \endcell% end previous cell contents
   \begingroup%
   \let\dg@getnodesize\tb@getnodesize% substitute routine to get nodesize
   \dg@USERSIZE=#1\relax%
   \ifnum\dg@USERSIZE<\@ne \dg@USERSIZE=\@ne \fi%
   \dg@parse{#2}%
   \dg@label{\tb@draw{#1}{#2}}}% draw arrow

\def\tb@getnodesize#1#2#3#4#5{\dimen3=\tb@arrowpad\dimen2 #4=\dimen3 #5=\dimen3\relax}
\def\tb@getnodesize#1#2#3#4#5{\ifnum#2=0\ifnum#3=0\tb@getnodesizetail{#4}{#5}\else\tb@getnodesizehead{#4}{#5}\fi\else\tb@getnodesizehead{#4}{#5}\fi}
\def\tb@getnodesizetail#1#2{\dimen3=.5\dimen2 #1=\dimen3 #2=\dimen3}
\def\tb@getnodesizehead#1#2{\dimen3=.5\dimen2 #1=\dimen3 #2=\dimen3}

\def\tb@draw#1#2#3#4{%
        \dg@X=0\dg@Y=0\dg@XGRID=1\dg@YGRID=1\unitlength=.001\dimen0%
        \dg@LBLOFF=\dgLABELOFFSET \divide\dg@LBLOFF\unitlength%
        \dg@drawcalc% compute arrow geometry
        \begincell% start tableau cell
        \let\lams@arrow\tb@lams@arrow% substitute routine
        \begin{picture}(0,0)\begingroup\dg@draw{#1}{#2}{#3}{#4}\end{picture}%
        \endcell% end tableau cell
        \endgroup% match \begingroup in \tb@arrow
        \begincell}% start new entry in this cell
}

\def\tb@lams@arrow#1#2{%
 \lams@firstx\z@\lams@firsty\z@
 \lams@lastx#1\relax\lams@lasty#2\relax
 \lams@center\z@
 \N@false\E@false\H@false\V@false
 \ifdim\lams@lastx>\z@\E@true\fi
 \ifdim\lams@lastx=\z@\V@true\fi
 \ifdim\lams@lasty>\z@\N@true\fi
 \ifdim\lams@lasty=\z@\H@true\fi
 \NESW@false
 \ifN@\ifE@\NESW@true\fi\else\ifE@\else\NESW@true\fi\fi
 \ifH@\else\ifV@\else
  \lams@slope
  \ifnum\lams@tani>\lams@tanii
   \lams@ht\ten@\p@\lams@wd\ten@\p@
   \multiply\lams@wd\lams@tanii\divide\lams@wd\lams@tani
  \else
   \lams@wd\ten@\p@\lams@ht\ten@\p@
   \divide\lams@ht\lams@tanii\multiply\lams@ht\lams@tani
  \fi
 \fi\fi
 \ifH@  \lams@harrow
 \else\ifV@ \lams@varrow
 \else \lams@darrow
 \fi\fi
}

\catcode`\@=\savecatcodeat
\let\savecatcodeat\undefined
\begin{document}    
\large

\title{On the denominators of Young's seminormal basis } 
\author{Steen Ryom-Hansen}

\author{
  Steen Ryom-Hansen{\thanks{Supported in part by FONDECYT grant{\color{black}s} 1171379  
and {\color{black}1221112 and by MathAmSud grant ALGonCOMB 22-MATH-01}}}}

%% \address{Instituto de Matem\'atica y F\'isica, Universidad de Talca \\
%% Chile\\ steen@@inst-mat.utalca.cl }
%% \footnote{Supported in part by FONDECYT grant 1090701.} 

%\date{}
\date{\vspace{-5ex}}
\maketitle

\begin{abstract}
{We study the seminormal basis $ \{ f_\T \} $ for 
the Specht modules of the Iwahori-Hecke algebra $ {\cal H}_{n}(q)$ 
of type $ A_{n-1} $. 
We focus on the base change coefficients 
between the seminormal basis $ \{ f_\T \} $ and Murphys' standard basis $ \{ x_\T \} $
with emphasis on the denominators of these coefficients.
In certain important cases we obtain simple formulas for these coefficients 
involving hook lengths.
Even for general {\color{black} standard} tableaux
we obtain new formulas.
On the way we prove a new result about submodules of the restricted Specht module at root of unity.}
\end{abstract}

\section{Introduction}

In this work we are concerned with the representation theory of the (Iwahori)-Hecke algebra 
of type $A$.

\medskip
As is well known, when defined over the fraction field
$ \cal K $ of ${\cal A}:= \Z[q, q^{-1}] $ 
the Hecke algebra $ {\cal H }^{\cal K}_n(q) $ 
is semisimple and its irreducible modules 
are the Specht modules $ S_q^{\cal K}({\lambda}) $
with $ \lambda $ varying over the set of 
integer partitions $ \Par_n$ of $ n$.
Each $ S_q^{\cal K}({\lambda}) $ is endowed with the
{\it seminormal basis} $ \{ f_\T \, | \, \T \in \std(\lambda) \}$, on which the action of the
Hecke algebra generators $ T_i $ can be described in terms of simple formulas,
{\it Young's seminormal form}, that have been known
for a long time.

\medskip
Our main interest is however rather the representation theory of
the Hecke algebras 
$ {\cal H }^{ {\Bbbk} }_n(q) $, defined over a field $ \Bbbk$ such that 
$ q $ is specialized to a root of unity in $   \Bbbk$. 
These specialized Hecke algebras are in general not semisimple, and
their representation theory is 
much more complicated than that of $ {\cal H }^{\cal K}_n(q) $, with many fundamental
problems still open. The
group algebra $ \F_p \Sy_n$ of
the symmetric group over the finite field $ \F_p$ 
is a special case of such an $ {\cal H }^{ {\Bbbk} }_n(q) $, and in spite of much progress in recent years
there is still no efficient algorithm for calculating the dimensions of the irreducible modules
for $ \F_p \Sy_n$.

\medskip
Let $ {\cal H }^{\cal A}_n(q) $ be the Hecke algebra defined over $ \cal A$. 
Then there are also Specht modules $  S_q^{\cal A}({\lambda}) $ for 
$ {\cal H }^{\cal A}_n(q) $
but for these {\it integral} Specht modules the seminormal basis does not exist and one
needs to work with Murphy's standard
basis, $ \{ x_\T \, | \, \T \in \std(\lambda) \} $, 
on which the action of $ T_i$ can only be described indirectly via a complicated recursion. On the
other hand, the standard basis has the advantage that it exists
for all specializations, including $ S_q^{\cal K}(\lambda) $. 

\medskip
In this work we
consider $ S_q^{\cal K}({\lambda}) $ and 
study the very natural question of determining the coefficients
%{\color{green}{\sout{between the $ \{ f_{\color{green}\T} \}  $-basis and
%      the $ \{ x_{\color{green}{\T}} \}  $-basis,}}}
{\color{black}{{of the $ \{ f_{\color{black}\T} \}  $-basis when 
expanded in the $ \{ x_{\color{black}{\T}} \}  $-basis,}}}
that is the base change matrix between the two bases.
For certain important {\color{black} standard} $ \lambda$-tableaux, that we call
{\it generalized James-Murphy tableaux}, we are able to get exact formulas for these coefficients.
The denominators that appear in the formulas involve certain hook lengths, and
when these denominators are nonzero in $ \Bbbk $, the corresponding seminormal basis
element will also exist in $ S_q^{\Bbbk}({\lambda}) $, which will have certain consequences for the
restriction of $ S_q^{\Bbbk}({\lambda}) $ from
$ {\cal H }^{\Bbbk}_{n}(q) $
to $ {\cal H }^{\Bbbk}_{n-1}(q) $, that we investigate in the final section of our paper.

\medskip
For completely 
general {\color{black} standard} $\lambda$-tableaux our methods do not give rise to a closed formula for
the base change matrix, but even in this case we obtain a 
fast algorithm for expanding the $ f_\T $'s in terms of
the $ x_\T$'s, that we explain.

\medskip
This paper has a long story. The first version of it was published on arXiv in 2009, and was formulated
only in the symmetric group setting, but in 2010
we uploaded a version of it to the arXiv in which 
our results were generalized to the Hecke algebra setup. 
The present version of the paper is essentially the same as the 2010 version, although we have corrected
a large number of errors and inaccuracies, and have added some new examples.

\medskip
The original motivations for our paper were twofold.
In \cite{RH} we showed that the coefficients of the quantum 
group action
%{\color{green}{\sout{of}}}
{\color{black}{on}} the Fock space, a main
ingredient in the LLT-algorithm for calculating decomposition numbers
for the Hecke $ {\cal H}^{\C}_{n}(q)$ with $ q \mapsto e^{2\pi i/l}$, see \cite{LLT} and \cite{Ar}, 
are closely related to the $ \{ f_\T \} $-basis. Indeed, 
let $ \T_n $ be the 
$ \lambda$-tableau that has $ n $ in a removable node and has the remaining 
numbers $ \{1,2, \ldots, n-1 \} $
filled in along rows; here is an example with $ n=26$
\begin{equation}\label{masexamplesTableaux_b}
      \T_{26} =  \raisebox{-.5\height}{\includegraphics[scale=0.7]{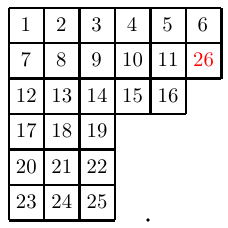}}
\end{equation}

This kind of $ \lambda $-tableaux plays an important role
in James and Murphy's calculation of the 
Gram matrix of the Specht module, 
\cite{JM}. Let $ \langle \cdot, \cdot \rangle $
be the canonical bilinear form on $  S_q^{\cal K}({\lambda}) $. We then 
proved in \cite{RH}, that 
the 
coefficients of the above mentioned quantum group action on the Fock space are the
$ l $-adic quantum valuation of the $ \langle f_{\T_n}, f_{\T_n} \rangle$'s. 
This observation made us speculate on a connection between the 
representation theory of $ {\cal H }_{n}^{\Bbbk}(q)  $
and the 
$ \{ f_{\color{black}\T} \} $-basis, and this was the first motivation for studying in more detail the expansion of
$ f_{\T_n}$.

\medskip
A second motivation for our work came from the strong analogy with the theory of symmetric functions.
The construction of the $ \{ f_{\color{black}{\T}} \} $-basis
is 
parallel to the construction of the Macdonald polynomials. 
Both are obtained through 
a Gram-Schmidt process
over a partial order which must first be extended to a total order
to perform the Gram-Schmidt process.
In the case of Macdonald polynomials {\color{black}{the initial basis is the basis of monomial
    symmetric functions}}, 
in the case of the seminormal basis the initial basis is the $ \{x_\T \} $-basis. 
By Cherednik's work, see \cite{C}, the 
Macdonald polynomials are independent of this extension because they are eigenvectors of 
operators coming from the double affine Hecke algebra; in the case of the seminormal basis 
this role is played by the Jucys-Murphy elements, 
see \cite{M3}. Finally, the formula for $ \langle f_{\T}, f_{\T} \rangle$ has
a striking similarity with its Macdonald polynomial analogue, 
see \cite{JM} and \cite{C}. 
But in the above picture an analogue for $ \{ f_\T \}$ of the 
positivity theory for the expansion coefficients of the 
Macdonald polynomials, due to Haiman and others, see \cite{H}, is 
still missing.
Our second motivation for the paper was to explain an attempt on how to fill this gap.

\medskip
Since the previous version of this paper, a number of papers treating related topics
have appeared, although the present paper is still the only one
which works at the Hecke algebra level of 
generality.
We here mention the paper by Raicu, \cite{Rai}, whose Theorem 1.2 is related to
our Theorem {\ref{f_n}}. 
Especially relevant are the two papers by
Fang, Lim and Tan \cite{FLT1} and \cite{FLT2}, in which the question of 
the splitting over $ \Z_p$ of the canonical map
$\Delta(\lambda + \mu)  \rightarrow  \Delta(\lambda) \otimes  \Delta(\mu) $ of Weyl modules for
the Schur algebra is treated. The authors connect this question to another question about
the denominators of Young's seminormal basis, and are this way able to answer it  
for certain pairs of partitions $\lambda, \mu$. 
It should be pointed out that the 
methods used in all the mentioned papers are different from ours, and different between them.

\medskip
{\color{black}{
We finally mention that computational evidence indicates that exact 
formulas for expressing the inverse coefficients of the  $ \{ x_\T \} $-basis in terms of
the $ \{ f_\T \} $-basis are more difficult to
get by. 
On the other hand, these inverse coefficients are the topic of a recent 
paper by Armon and Halverson, \cite{AH}, who use results of Ram
to give a fast algorithm for calculating
them.}}

\medskip
Let us now explain in more detail the contents of the paper.
In section \ref{Basic notations and results} we give a precise
formulation of the setting in which we shall be working,
and recall the relevant results from the literature, the most important being
Murphy's standard basis
$ \{ x_\T \, | \, \T \in \std(\lambda) \} $ for the Specht module 
with associated Hecke algebra action
{\eqref{action}} and the Garnir relations {\eqref{Garnir}}. 
We also recall Young's seminormal basis
$ \{ f_\T \, | \,  \T \in \std(\lambda) \} $ 
and
its Hecke algebra action, that is Young's seminormal form, or YSF, see {\eqref{Seminormal_form}}.
In section \ref{analgorithm} we explain a simple algorithm for calculating
the expansion of $ f_\T $ in terms of $ x_\s$'s, but where $ \s $ runs over
all {\color{black}row} standard tableaux, not just standard tableaux. 
In section \ref{A generalization of the James-Murphy tableaux} we introduce the generalized James-Murphy
tableaux that enter in our main Theorems. In section \ref{YSF along one row}
we describe a formula that results from applying YSF along one row of the generalized James-Murphy tableau.
In section \ref{Fat hook partitions} we treat the expansion
of $ f_\T $ for $ \T $ a fat hook tableaux. The
formula of the previous section \ref{YSF along one row}
has a certain similarity with  
the Garnir relations and this similarity is a main ingredient of Theorem 
\ref{fat_hook_lemma} of that section.
In section \ref{relevantsection} we treat the expansion of the
seminormal basis elements corresponding to all generalized
James-Murphy tableaux, in essence, by reducing to the case of fat hook partitions. 
In the final section \ref{restricted} we consider the expansion in the case of
general {\color{black}standard} tableaux.
We also give an application of the theory developed in the previous sections to
the modular Hecke algebra representation theory, that is the representation theory of
$  {\cal H }^{ {\Bbbk} }_n(q) $, giving a criterion for a certain Specht module to split off from
the restricted Specht module $ \Res \, S_q^{\Bbbk}({\lambda}) $. This criterion was also present in
the previous versions of our article, but the statement and proof of it were very inadequate. 
The present formulation of the criterion is inspired by the formulation of the splitting criterion in \cite{FLT1}, 
but our proof technique is quite different from the one used in \cite{FLT1}.

\medskip
    {\color{black}
  It is a great pleasure to thank M. Fang, K. J. Lim and K. M. Tan for useful email conversations related to
  this article.
It is a special pleasure to thank the referee
for suggesting Corollary \ref{cor} and for providing us 
with a detailed list
of comments that helped us
eliminate errors and improve the presentation.
}

\section{Basic notations and results}\label{Basic notations and results}
In this section we introduce the notation and recall some basic results that shall be used throughout the paper.

\medskip
Let ${\cal A}:= \Z[q,q^{-1} ] $ be the ring of Laurent polynomials over $ \Z$
and let $ \cal K $ be its quotient field. 

\medskip
We denote by $ {\cal H }^{\cal A}_{n}(q)  $ the {\it Iwahori-Hecke algebra} of type $ A_{n-1} $ over $ \cal A$.
That is, 
$ {\cal H }^{\cal A}_{n}(q)  $ is the $ \cal A $-algebra on generators $ T_1, T_2,  \ldots, T_{n-1} $ 
subject to the relations 
$$ \begin{array}{ll}
T_i T_j = T_j T_i & \mbox{for } \, | i-j | > 1 \\ 
T_i T_j T_i = T_j T_i T_j & \mbox{for } \, | i-j | = 1 \\
(T_i -q ) (T_i +1 ) = 0 & \mbox{for } \, i= 1,2, \ldots , n-1.
\end{array}
$$
For an $ \cal A $-algebra $ {\cal B} $ we introduce the {\it specialized Iwahori-Hecke} algebra 
$ {\cal H }^{{\cal B}}_{ n}(q) :=     { \mathcal B }\otimes_{\cal A} {\cal H }^{\cal A}_{n}(q)   $.
If $ {\cal B} = \cal K $ we also sometimes write $ {\cal H }_n = {\cal H }^{\cal K}_{ n}(q) $.
In general, 
we shall refer to $ {\cal H }^{\cal A}_{n}(q)  $ and to the
various variations of it simply as the {\it Hecke algebra},
omitting the name Iwahori.

\medskip
Let $ \Sy_n $ be the symmetric group on $ n $ letters.
{\color{black}{There is a natural action of $ \Sy_n $ on $ \{1, 2, \ldots, n \} $ which
we view as a right action}}. 
It is well known that  $ \Sy_n $ is a Coxeter group on 
the set $ S:=\{ s_1, s_2, \ldots, s_{n-1} \} $ where $ s_i$ is the simple
transposition $ s_i := (i,i+1 ) $.
For $ w = s_{i_1} s_{i_2} \cdots s_{i_N} \in \Sy_n    $ a reduced expression for $ w $ we set
$ T_w := T_{i_1} T_{i_2} \cdots T_{i_N} \in {\cal H }^{\cal A}_{n}(q) $. Then it follows from Matsumoto's Theorem
that $ T_w $ is independent of the choice of reduced expression for $ w$. Moreover, as is also well known, 
$ \{ T_w \, | \, w \in \Sy_n \} $ is an $ \cal A $-basis for $ {\cal H }^{\cal A}_{n}(q) $
and $ \{1 \otimes  T_w \, | \, w \in \Sy_n \} $ is a $ \cal B $-basis for
$ {\cal H }^{\cal B}_{n}(q) $.

\medskip
Any field $ {\Bbbk} $ with identity $ 1 $ can be made into an $ \cal A$-algebra via
$ q \mapsto 1 $ and in this case we have that ${\cal H }^{\Bbbk}_n(q) \cong \Bbbk \Sy_n$, where 
$ {\Bbbk} \Sy_n $ is the group algebra of the symmetric group. Thus all results that hold for specialized
Hecke algebras also hold for $ {\Bbbk} \Sy_n $.

\medskip
We next recall the combinatorial notions of partitions and tableaux, associated with the Hecke algebra.
For $ n $ a positive integer we  
denote by $ \Par_n $ the set of {\it integer partitions} of $ n $. 
To be precise, an element $ \lambda $ of $\Par_n $ is a weakly decreasing sequence of positive integers
$ \lambda = (\lambda_1 , \lambda_2, \ldots,  \lambda_K ) $ such that
$ \lambda_1 +\lambda_2 + \ldots + \lambda_K = n$. Similarly, we denote by $ \Comp_n $
the set of {\it integer compositions} of $ n $, defined the same way as $\Par_n$ but
allowing each $ \lambda_i $ to be zero and
omitting the condition that $ \lambda $ be weakly decreasing.
We set $ \Par:= \cup_n \Par_n $ and $ \Comp:= \cup_n \Comp_n $. 

\medskip
The sets $ \Par_n $ and $ \Comp_n $ are endowed with order relations that play
an important role in the representation theory of the Hecke algebra, and also elsewhere.
For $ \lambda, \mu \in
\Comp_n $ we may assume that there is a $ K $ such that 
$ \lambda = (\lambda_1 , \ldots, \lambda_K ) $ and $ \mu = (\mu_1 ,  \ldots , \mu_K )  $, by extending
with zeros if necessary. Then the {\it dominance order} on $ \Comp_n $ is defined via 
$ \lambda \unlhd \mu $ if 
\begin{equation} \sum_{i=1}^j \lambda_i  \leq \sum_{i=1}^j \mu_i, \, \, \, \,
  \mbox{ for all }  1 \le j \le K. \end{equation}
By restriction, it induces an order relation on $ \Par_n$, denoted the same way.

\medskip
The dominance order is only a partial order on $ \Par_n $ and $ \Comp_n $ but can be 
extended to a total order on both sets, for example via
$ \lambda < \mu $ if there is a $ j $ such that 
$ \sum_{i=1}^j \lambda_i  <   \sum_{i=1}^j \mu_i  $ but 
$ \sum_{i=1}^{j^{\prime}} \lambda_i  =	    \sum_{i=1}^{j^{\prime}} \mu_i  $ for all
$ j^{\prime} < j $.

\medskip
Let $ \lambda = (\lambda_1,\ldots,  \lambda_K ) \in \Par_n $.
Then the {\it Young diagram} $ {\cal Y }(\lambda) $ 
for $ \lambda $ is the graphical representation of
$ \lambda $ through boxes, called {\it nodes}, arranged in $ K $ left adjusted lines,
with $ \lambda_1 $ nodes in the first line, $ \lambda_2 $ nodes in the second line just below the
first line, and so on.
For example, if $ \lambda = (5,5,4,1,1) \in \Par_{\color{black}{16}}$ and
$ \mu = (2,0, 3, 4,1) \in \Comp_{\color{black}{10}} $ we have
that
 
\begin{equation}\label{Youngex}
{\mathcal Y(\lambda)} =   \raisebox{-.5\height}{\includegraphics[scale=0.7]{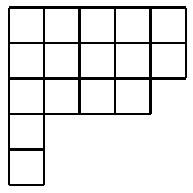}}, 
\, \, \, \, \, \, \, \, \, \, \, \, \, \, \, \,
{\mathcal Y(\mu)} =   \raisebox{-.5\height}{\includegraphics[scale=0.7]{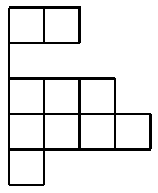}}
\end{equation}
For $ \lambda \in \Comp_n $, and in particular for $ \lambda \in \Par_n $,
we use matrix convention when referring to the nodes of 
$ {\cal Y }(\lambda) $. In other words $ [r, c ] $ refers to the node that occurs in $ {\cal Y }(\lambda) $
in the $ r$'th row from the top and the $ c$'th column from the left. A 
$ \lambda $-tableau $ \T $ is a filling of the nodes of $ {\cal Y }(\lambda) $ with 
the numbers $ 1, 2, \ldots, n$, each number occurring exactly once. 
We write $ \T[r,c] $ for the number that occurs in the $ [r,c] $'th node of $ \T $.
Formally, $\T $ is a bijection from the nodes of $ {\cal Y }(\lambda) $ to the numbers
$ \{ 1, 2, \ldots, n \} $ and so $ \T^{-1}(i) $ is the node of $ {\cal Y }(\lambda) $ that contains $ i$.
If $ \T $ is a $ \lambda $-tableau we say that $ \T $ is {\it row standard} 
if the numbers of each row of $ \T $ appear increasingly from left to right, and 
that $ \T $ is {\it column standard} if the numbers of each column of $ \T $ appear increasingly from top 
to bottom. If $\T$ is both row and column standard, we say that is {\it standard}.
We let $ \T \mapsto shape(\T) $ be the function that 
maps $ \T $ to its underlying partition, thus $ shape(\T)= \lambda $ if $ \T $ is
a $ \lambda$-{\color{black}tableau}. We denote by $ \tab(\lambda) $ and $ \std(\lambda) $ the sets of
$ \lambda$-tableaux, and 
standard
$  \lambda$-tableaux, respectively, and similarly we set $ \std(n):= \cup_{\lambda\in \Par_n} \std(\lambda) $
and $ \tab(n):= \cup_{\lambda\in \Par_n} \tab(\lambda) $.

\medskip
For $ \T \in \std(\lambda)$ and $ 1 \le k  \le n $ we let $ \T|_{ 1,2, \ldots, k }  $ denote 
the tableau obtained by deleting the nodes of ${\cal Y}( \lambda)$
containing the numbers $ k+1, k+2, \ldots, n $. With 
this notation the dominance order on $ \Par_n $ 
is extended to {\color{black}standard} $ \lambda$-tableaux as follows: $ \s {\color{black}\unlhd} \T $ if
$ shape(\s|_{ 1,2, \ldots, k  }) \leq  shape(\T|_{ 1,2, \ldots, k }) $ for all $ k $. 
In a similar way, the total order $ < $ is extended to a total order on 
$ \std(\lambda)$, that is $ \s < \T $ if $ shape(\s|_{ 1,2, \ldots, k  }) <  shape(\T|_{ 1,2, \ldots, k }) $ for some $ k $
and $ shape(\s|_{ 1,2, \ldots, k  }) =  shape(\T|_{ 1,2, \ldots, k }) $ for all $ k^{\prime} < k $.
There is a unique maximal $ \lambda $-tableau $ \T^{\lambda} $ with respect to both orders;
it has the numbers $ 1,2, \ldots, n $ filled in increasingly along 
the rows. There is also a unique minimal $ \lambda $-tableau $ \T_{\lambda} $
with respect to both orders; it 
has the numbers $ 1,2, \ldots, n $
filled in increasing{{\color{black}ly} along the columns. Here are examples of $ \T^{\lambda} $ and $ \T_{\lambda} $,
where $ \lambda $ is as in $\eqref{Youngex}$.

\begin{equation}\label{tabex}
\T^{\lambda} =   \raisebox{-.5\height}{\includegraphics[scale=0.7]{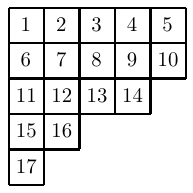}}\, \, \, , 
\, \, \, \, \, \, \, \, \, \, \, \, \, \, \, \,
\T_{\lambda}=   \raisebox{-.5\height}{\includegraphics[scale=0.7]{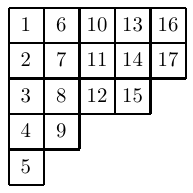}}
\end{equation}

\medskip
Since $  \Sy_n $ is a Coxeter group it is endowed with the Chevalley-Bruhat order $ <$.
We choose the convention that $ 1 \in \Sy_n $ is the maximal element with
respect to $ <$. The following compatibility between the Chevalley-Bruhat order and the dominance order 
is known as Ehresmann's Theorem, see \cite{Ma} for a proof. 
\begin{Thm}\label{ehresmann}
  Let $ \lambda \in \Par_n$ and suppose that $ \s, \T $ are
  standard $ \lambda$-tableaux.
  Then $ \s \lhd \T $ if and only if $ d(\s) < d(\T)$. 
\end{Thm}  

\medskip
We now recall the construction of Murphy's standard basis for the Hecke algebra.
$ \Sy_n $ acts {\color{black}on the right}
on the 
set of $ \lambda$-tableaux $ \T $ by permuting the numbers inside each $ \T $,
%% {\color{green}{\sout{
%% where we adopt the convention that this
%% action is a {\it right action}}}}
{\color{black}{that is, formally, $ (\T w)([r, c ]):= \T ([r, c ])w$ for $ \T \in \tab(\lambda) $ and
$w \in \Sy_n$}.}
For $ \T $ a $ \lambda$-tableau, we introduce $d(\T) \in \Sy_n $ by requiring that 
$ \T = \T^{\lambda} d(\T) $ and let $ \Sy_{\lambda} \subseteq \Sy_n $ denote the row 
stabilizer of $ \T^{\lambda} $
under the $ \Sy_n$-action. 
For a pair $  (\s, \T) $ of $ \lambda$-tableaux, 
Murphy introduced in this context the following elements of $ {\cal H }^{\cal A}_n(q)$
\begin{equation}\label{Murphy's theory} x_{\lambda} := \sum_{w \in \Sy_{\lambda} } T_w, \, \, \, \, \, \, \, \, 
  x_{\s \T} := T_{d(\s)^{-1} }\, x_{\lambda} \,T_{d(\T)}. 
\end{equation}
If $ \s $ and $ \T $ are standard tableaux we say that $   x_{\s \T} $ is a {\it standard element}, otherwise
we say that it is a {\it nonstandard element}.
Murphy proved in \cite{M1} that the standard elements $   x_{\s \T} $
form an $\cal A $-basis for $ {\cal H }^{\cal A}_{n}(q) $, the {\it standard basis}. 
They also induce a basis for the specialized Hecke algebra.

\medskip Let $ \lambda \in \Par_n$ and let
$ \overline{N}_{\lambda} $ be the $ \cal A $-span of $\{ x_{\s \T} \, | \, \s, \T \, \mbox{ are } 
\mu \mbox{-tableaux}\, \mbox{ with } 
\mu \rhd \lambda \} $. 
Then Murphy showed that $ \overline{N}_{\lambda} $ is a two-sided ideal of $ { \cal H}^{\cal A}_n(q) $ and 
defined the  
{\it Specht module} $ S_q^{\cal A}({\lambda}) $ for 
$ {\cal H}^{\cal A}_n(q) $ as the right
submodule of $ {\cal H}^{\cal A}_n(q)/ \overline{N}_{\lambda}$, generated by 
$ x_{\lambda} + \overline{N}_{\lambda} $.
It is a generalization of the {\it Specht module} known from the representation
theory of $ \Sy_n$, or more precisely of the {\it dual Specht module}.

\medskip
    {\color{black}In the representation theory of $ \Sy_n$,
      as exposed for example in \cite{J},
      the Specht modules play a key role. For a field $ \Bbbk $ that contains $ \Q$,
      the Specht modules are the simple modules 
      for $ \Bbbk \Sy_n$, and even for arbitrary fields $ \Bbbk $ they
      can be used to describe the simple modules for $ \Bbbk \Sy_n$, at least in principle.
      Indeed, the Specht modules are endowed with
      certain, combinatorially defined, bilinear forms, and each simple module for $ \Bbbk \Sy_n $
      can be realized as a quotient of a unique Specht module by the radical of its bilinear form.
      Similar statements also hold for the dual Specht modules for $ \Sy_n $, but it should
      be noted that the parametrizations of the simple modules for $ \Bbbk \Sy_n $, using 
    the Specht modules or the dual Specht modules, are different. 

\medskip   
Returning to the Specht module $  S_q^{\cal A}({\lambda}) $ for
${\cal H}^{\cal A}_n(q) $, we have that}
it is free over $ \cal A $ with basis 
given by $  x_{ \T } := x_{ \T^{\lambda} \T } + \overline{N}_{\lambda} $
where $ \T $ runs over $ \std(\lambda) $. 
We refer to this basis as the {\it standard basis} for $ S^{\cal A}_q({\lambda}) $. 
For $ \cal B $ 
an $ \cal A $-algebra, the standard basis induces a basis for the specialized Specht module
$ S_q^{\cal B}({\lambda}) : = {\cal B } \otimes_{\cal A} S^{\cal A}_q(\lambda) $.
If $ {\cal B } = \cal K $ we shall also sometimes write 
$ S({\lambda}) = S_q^{\cal K}({\lambda})$.

\medskip
The standard basis $ \{ x_{\s \T} \}$ is a 
{\it cellular basis} for $ {\cal H }^{\cal A}_n(q) $ in the sense of Graham and Lehrer, see \cite{GL} and \cite{Ma}.
Thus $ S^{\cal A}_q({\lambda}) $ is endowed with a 
symmetric bilinear form $ \langle \cdot, \cdot \rangle_{\lambda} $. It is concretely 
given by 
$ \langle x_\s, x_\T \rangle_{\lambda}  = a $ where $ a $ is the coefficient of $ x_{\lambda } $ in 
$ x_{\lambda} T_{d(\s)} T_{d(\T)^{-1}} x_{\lambda}  $
when expanded in the standard basis.	
For $ \cal B = K $ the form $ \langle \cdot, \cdot \rangle_{\lambda} $
is nondegenerate, 
but for example for $ {\cal B} = \F_p  $, the finite field of characteristic $ p $, 
the form may be singular.
In any case, the radical $  \rad \langle \cdot  , \cdot  \rangle_{\lambda} $ is 
a submodule of $ S^{\cal B }_q({\lambda}) $. {\color{black}For $ {\mathcal B} = \Bbbk $ a field,
we have that 
$ S^{\Bbbk }_q({\lambda}) / \rad \langle \cdot  , \cdot  \rangle_{\lambda} $ is either simple or zero and 
the nonzero modules that arise this way provide a classification of the simple modules for 
$ {\cal H}^{\Bbbk}_n(q) $.}		

\medskip
We now describe the action of the $ T_i $'s on the standard basis for $ S^{\cal A}_q({\lambda}) $.
Assume that $ \T \in  \std(\lambda)$ 
and that $ \s = \T s_i$. 
The action of $ {\cal H}^{\cal A}_{n}(q) $ on $ x_{ \T} $ is then 
given by the following formulas 
\begin{equation}{\label{action}}
  x_{ \T} T_i := \left\{ 
\begin{array}{ll} 
q x_{ \T} & \mbox{ if }  i \mbox{ and } i+1 \mbox{ are in the same row of } \T  \\
  x_{ \s } & \mbox{ if }  i +1 \mbox { is in a row below } i \mbox{ in } \T   \\
q x_{ \s} + (q-1)x_{ \T} & \mbox{ if } i +1 \mbox { is in a row above } i \mbox{ in } \T.
\end{array} \right.
\end{equation}
Unfortunately, $ \s $ may be nonstandard even if $ \T $ is a standard 
and so we  
need straightening rules to express nonstandard elements in terms
%{\color{green}{\sout{of}}}
of standard
elements.

\medskip
The relevant straightening rules are the $q$-analogues of the {\it Garnir relations}, generalizing the
usual Garnir relations known 
from the representation theory of the symmetric group $\Sy_n$. Let $ \lambda \in \Par_n $, and 
choose
$ (i,j) $ such that 
$ i \geq 1 $ and $ j \leq \lambda_{i+1} $. Define $ \mu := (\lambda_1, \ldots, \lambda_{i-1}, j-1, j ) $
{\color{black}{if $ i > 1 $ and $ \mu := ( j-1, j ) $ if $ i = 1 $ }}
    and suppose that
$ \mu \in \Comp_m $.
    Then the $ (i ,j) $-{\it Garnir}
    tableau $ \g_{i     j} $ is the $ \lambda $-tableau
satisfying that $  \g_{i  j}|_{[1,2, \ldots, m ]} = \T^{\mu} $ and 
that the numbers $ m+1, m+2, \ldots, n $ are filled in increasingly along the rows in the difference 
$ {\cal Y}(\lambda) \setminus {\cal Y}(\mu) $. In particular, $ \g_{{i j}} $
  is not column standard, 
since there is a descent between the nodes $ {\color{black}{[}}i, j
  {\color{black}{]}} $ and
$ {\color{black}{[}}i+1, j
  {\color{black}{]}} $.
Below we give the example $ \g_{{\color{black}2},3} $ using the same partition $\lambda=(5,5,4,2,1) $ as before.

\begin{equation}\label{garnir}
\g_{{\color{black}{2}},3} =   \raisebox{-.5\height}{\includegraphics[scale=0.7]{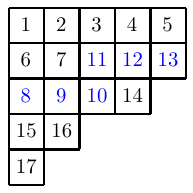}}\, \, \, 
\end{equation}

Set $ k:= \g_{{\color{black}ij}}[i{\color{black}{+1}},1] $. Then the numbers $ \{ k, k+1, \ldots, m \} $ are the numbers of the
$(i,j)$-{\it Garnir belt}, in \eqref{garnir} they are 
coloured blue. Let $ \Sy_{k, m } $ be the subgroup of 
$ \Sy_n $ consisting of the elements that fix 
$ \{ 1, 2 \ldots , n \} \setminus \{ k,k+1, \ldots, m \} $ pointwise. 
Set
\begin{equation}
  \G_{ij} :=
  \left\{ \g_{ij} w \in \tab(\lambda) \, \mid \, w \in \Sy_{k, m }\, \mbox{and} \, \g_{ij} w \mbox{ is row-standard}
    \right\}.
\end{equation}  
The $q$-analogue of the {\it Garnir relation} is now the following relation in $ S^{\cal A}_q({\lambda})  $
\begin{equation}{\label{Garnir}}
  \sum_{ \s \in \G_{ij}} x_{\s}     = 0. 
\end{equation}
The only nonstandard tableau appearing in $  \G_{ij}$ is $ \g_{ij} $ itself, and so {\eqref{Garnir}} can be used
to express $ x_{\g_{ij}} $ in terms of standard elements.
Using arguments explained in \cite{M2}, this can be extended to all nonstandard
$ x_\T$.

\medskip
We next explain the {\it Jucys-Murphy elements} for $ {\cal H }^{\cal A}_n(q)$.
For $ m= 1,2,3, \ldots, n $ we define 
$ L_m \in {\cal H}^{\cal A}_{n}(q) $ via 
\begin{equation}
 L_m := q^{-1}\, T_{(m-1, m)} + q^{-2} \,T_{(m-2, m)} + \ldots + q^{1-m}\,  T_{(1, m)}
\end{equation}
where $ (i,j )  $ is the element of $ \Sy_n $ given in usual cycle notation,
and where by convention we set
$ L_1 :=0 $. These are the Jucys-Murphy elements for $ {\cal H}^{\cal A}_{n}(q) $. For $ k \in \Z $ we introduce
the following {\it Gaussian integer} 
\begin{equation}
  [k]_q :=  \left\{\begin{array}{ll} 1+ q +q^2 + \cdots + q^{k-1} & \mbox{ if } k>   0 \\
{\color{black}0} & \mbox{ if } k = 0 \\
  -( q^{-1} +q^{-2} + \cdots + q^{-k}) & \mbox{ if } k <   0.  \end{array} \right.
\end{equation}
Thus for
%{{\color{green}{\sout{$ k \neq 0 $}}}
{\color{black}{all $k \in \mathbb Z $}}
we have $  [k]_q=  \frac{q^k -1}{q-1}  $ and 
{\color{black}$ \lim_{q \rightarrow 1}  [k]_q = k $}. 
With this at hand we have for 
$  \T \in \tab(\lambda) $ the {\it content function} 
\begin{equation}\label{content}
  c_{\T}: \{1,2,\ldots,n\} \rightarrow {\cal A}, \,  i \mapsto  [c-r]_q  \mbox{ where } \T[r,c] = i. 
\end{equation}
If $\T = \T^{\lambda} $ we often write $ \lambda$ for $ \T^{\lambda}$ as a subscript, for example
$ c_{\T}(i)= c_{\lambda}(i) $. 
In the context of the specialized Hecke algebra $ {\cal H}_n^{\cal B}(q) $, the content function 
$   c_{\T} $ is given by the same formula \eqref{content} but now takes values in $ \cal B$.
For $ \lambda \in \std(\lambda) $,
the $ L_i$'s satisfy the following JM-triangularity property with respect to the
$c_{\T}$'s. This key property is due to Murphy, see \cite{M3}, to our
knowledge it does not appear 
anywhere before \cite{M3}. 
\begin{equation}\label{key}
L_i x_\T = c_\T{\color{black}{(i)}} x_\T + \sum_{\s \in \std(\lambda),  \s  \rhd \T}   a_{\s} x_\s, \, \,   \mbox{ where } a_\s \in \cal {\cal A}.
\end{equation}  

\medskip
Let $ \cal B = \cal K$. Then 
an important application of the $ L_i$'s is the construction of idempotents in the Hecke algebra
$ {\cal H}_n $.
Let $ \Cont_n := \{   c_{\T}{\color{black}{(i)}} \, | \, \lambda \in \Par_n, \T \in \std(\lambda) \}$
and define for a $ \lambda$-tableau $ \T $ the following element of
$  {\cal H}_n $
\begin{equation}{\label{E-expansion}}
E_\T := \prod_{m=1}^n \prod_{c \in \Cont_n \setminus c_\T(m) }
\, \frac{L_m -c }{ c_\T(m) -c}.
\end{equation}
Then the $ E_{\T}$'s
are a complete set of primitive idempotents in ${ \cal H}_n$, 
for $ \T $ running over all standard tableaux. 
In case $ \T $ is a nonstandard tableau, the formula {\eqref{E-expansion}} also defines 
an idempotent $E_\T$ in $ {\cal H}_n$, but one gets nothing new since
$E_\T =  E_\s $ for some standard tableau $ \s $.
For $\T \in \std(\lambda) $, the idempotent $ E_\T $ gives rise to an
element $ f_{\T} $ of $ S_q^{\cal K}(\lambda)$ as follows
\begin{equation}
f_\T = x_\T E_\T \in  S_q^{\cal K}(\lambda). \, \, \, \,  
\end{equation}
Then $ \{ f_\T \, | \, \T \in \std(\lambda \} $ 
is the $q$-analogue of the {\it seminormal basis} for $ S_q^{\cal K}({\lambda}) $. 
It consists of common eigenvectors for the $ L_i$'s, with eigenvalues given by the contents, that is 
\begin{equation}{\label{dominant_terms}}
L_i f_\T = c_\T(i) f_\T, \, \, \,   i= 1, \ldots, n.  
\end{equation}
We keep $ {\cal B  } = \cal K $ and so we have that $ \langle \cdot, \cdot \rangle_{\lambda} $ is nondegenerate. 
Moreover the $ L_i $'s are selfadjoint with respect to $ \langle  \cdot, \cdot \rangle_{\lambda} $ 
and for any $ \s, \T \in \std(\lambda), \s \neq \T $, there is an $ i $ such that $ c_\s(i) \neq c_\T(i) $.
From this we deduce that the $ f_\T$'s are orthogonal with respect to $ \langle  \cdot, \cdot \rangle_{\lambda} $.
Furthermore we have the following triangular expansion, which is a consequence of \eqref{key}
\begin{equation}\label{ftriangular}
f_\T = x_\T + \sum_{\s \in \std(\lambda),  \s  \rhd \T}   a_{\s} x_\s, \, \,   \mbox{ where } a_\s \in \cal K.
\end{equation}

\medskip
The seminormal basis $ \{f_\T \} $ for $ S_q^{\cal K}(\lambda) $
can also be 
constructed using a Gram-Schmidt orthogonalization process on the standard basis $ \{x_\T \} $,
with respect to $ \langle \cdot, \cdot \rangle_{\lambda} $. 
Recall the extension of $ \lhd $ to a total order $ < $ on $ \std(\lambda) $.
For the Gram-Schmidt 
process we first take $ f_{{\lambda}} =  f_{\T^{\lambda}} := x_{{\lambda}} = x_{\T^{\lambda}} $
and then continue recursively downwards
along $ < $ as follows
\begin{equation}\label{downwards}
  f_{\T} := x_{\T} - \sum_{ \substack{\s \in \std(\lambda) \\  \s > \T}}
 \frac{\langle f_{\s}, x_{\T} \rangle_{\lambda}}{\langle f_{\s}, f_{\s} \rangle_{\lambda}}
 f_{\s}.
\end{equation}

Apriori, the orthogonal basis that results from this Gram-Schmidt process may depend on how
$ \lhd  $ is extended to a total order. 
On the other hand, using \eqref{key}
one checks that replacing $ <$ by $ \lhd$ in 
\eqref{downwards} one obtains the same $ f_\T$'s from the Gram-Schmidt process. 
In other words, in \eqref{downwards} we can use any extension of $ \lhd $ to a total order, without
changing the outcome.

\medskip
This formalism is reminiscent of basic principles in the theory of symmetric functions.
A natural basis 
for the space of symmetric functions $ \Sym$ is given by the monomial symmetric functions
$ \{  m_{\lambda} |  \lambda \in \Par  \} $. 
The Macdonald polynomials $ \{  P_{\lambda}  |  \lambda \in \Par  \} $
are another basis for $ \Sym$ which is 
constructed via 
a Gram-Schmidt algorithm on $ \{  m_{\lambda}  |  \lambda \in \Par  \} $, 
using a certain inner product on 
$ \Sym$. The order relation {\color{black}{here}} is an extension
of the dominance to a total order on all of $ \Par$.
On the other hand the $P_{\lambda}$'s 
can also be realized as common eigenvectors for a family of selfadjoint operators 
on $ \Sym $ that have their origin in the Cherednik algebra. As was the case for the seminormal basis, 
one then concludes that the Gram-Schmidt construction of $P_{\lambda}$ does not depend on the choice of
extension of the dominance order to a total order on $ \Sym$.

\medskip
It should be pointed out, however,
that in both cases it is difficult to gain information about the
orthogonal basis directly from the Gram-Schmidt process, and
so the Gram-Schmidt process 
is more a tool for calculating examples than a tool for deducing theoretical properties of the basis.

\medskip
Returning to the seminormal basis $ \{ f_{\T} | \T \in \std(\lambda) \}$, 
the action of the $ T_i $'s on $ S_q^{\cal K}(\lambda ) $ is much easier to describe
%{\color{green}{\sout{using it}}}
than using the standard basis, since the Garnir relations are not needed.
In fact this is one of the big advantages of the seminormal basis over the standard basis.

\medskip
Suppose that $ \T \in \std(\lambda) $ and let $ \s:= \T s_i $.  Define the {\it radial distance} 
$ \rho := c_{\T}(i) - c_{\s}(i) =  c_{\T}(i) - c_{\T}(i+1) $. Then we have the following formulas,
known as {\it Young's seminormal form}, or simply YSF.
They play an important role throughout the paper.

\begin{equation}{\label{Seminormal_form}}
f_\T \,T_{i} = 
\color{black}{\left\{\begin{array}{ll} 
  q f_\T & \mbox{ if } i  \mbox{ and } i+1  \mbox{ are in the same row of } \T  \\
  -f_\T & \mbox{ if } i  \mbox{ and } i+1  \mbox{ are in the same column of } \T  \\
    - \dfrac{1}{[\rho]_q} f_\T + f_\s   & \mbox{ if } \s  \mbox{ is standard and } i+1 \mbox{ is in a row below } i
\mbox{ in } \T
  \\
\dfrac{-1}{[\rho]_q} f_\T + 
\dfrac{q [\rho +1]_q [\rho {\color{black}{\, -1}}]_q}{[\rho ]_q^2} f_\s    & \mbox{ if } \s  \mbox{ is standard and }
 i+1 \mbox{ is in a row above } i
\mbox{ in } \T.
\end{array}
\right. }
\end{equation}
These formulas were proved in Theorem 6.4 of \cite{M1}, although the formulation there contained
several minor errors, see the discussion in the proof of Theorem 2.3 of \cite{FLT1}.

\section{An algorithm using Young's seminormal form}\label{analgorithm}

As mentioned in the introduction, our 
goal is to find a 
formula for the base change matrix between the $ \{f_{\T}\}$ basis and the $ \{ x_\T \}$
basis for $ S_q^{\cal K}(\lambda)$, or equivalently
to determine the $ a_\s$'s that appear in
\eqref{ftriangular}.

\medskip
A first observation is that Young's seminormal form,
that is formula {\eqref{Seminormal_form}}, in fact does give rise to
an algorithm for writing $ f_\T $ as a linear combination of $ x_\s $'s, but with  
$ \s $ running over {\it all of} $ \tab(\lambda) $, not just $ \std(\lambda) $.
In this section we explain this algorithm.

\medskip
Suppose that $ \T \in \std(\lambda)$ and write 
$ d(\T) = s_{i_1} s_{i_2} \cdots s_{i_k} $ in reduced form. From this we get a chain
of standard tableaux $ \{ \T_0, \T_1, \ldots, \T_k \} \subseteq \std(\lambda) $ as follows 
\begin{equation}{\label{tableauxchoise}}
\T_0 := \T^{\lambda}, \, \T_1 := \T_0 \, s_{i_1}, \, \T_2 := \T_1 \,s_{i_2}, \, \ldots, \, 
\T_k := \T_{k-1} \,s_{i_k}
\end{equation}
where $ \T = \T_k$. Using Theorem \ref{ehresmann} we get that $ \T_{j}  \lhd \T_{j-1} $ for all $j $,
and so $ i_{j}+1 $ appears in a row below $ i_j $ in $ \T_{j-1}$ for all $ j$. This means that
when calculating $ f_{\T_{j-1}}  T_{j} $ using YSF, we are in the third case of
\eqref{Seminormal_form}.

\medskip
For our algorithm we start out with the equality $ f_\lambda = x_\lambda$. 
Applying $ T_{i_1} $ to each side of this equality, the left hand becomes $ f_\lambda T_{i_1} $ which can
be calculated using YSF, whereas the right hand side
becomes $ x_\lambda T_{i_1} = x_{\T_1}  $, via {\eqref{action}}. To be precise, 
when applying YSF we get 
\begin{equation}
  \begin{array}{l}
  f_{\T_0} T_{i_1} =  x_{\T_1}  \Longleftrightarrow  f_{\T_0} + a_{1} f_{\T_1} = x_{\T_1}
  \Longleftrightarrow  f_{\T_1 }=   a_{\T_1}^{-1}(-  f_{\T_0} + x_{\T_1} )    \Longleftrightarrow \\ [.3cm]
  f_{\T_1 } =   a_{2} x_{\T_0} +  a_{3} x_{\T_1}  
  \end{array}
  \end{equation}
for certain explicit coefficients $ a_{1}, a_{2}, a_{3} \in \cal K^{\times}$. 
This is the first step of the algorithm. For the next step 
we apply $ T_{i_2} $ to both sides of the equation $ f_{\T_1 } =   a_{2   } x_{\T_0} +  a_{   3} x_{\T_1}   $, 
that we just found. The left hand side is $ f_{\T_1 } T_{i_2}$ that can be calculated using YSF and
the right hand can be calculated via {\eqref{action}}. In more detail, for the left hand side we
get via YSF that $ f_{\T_1 } T_{i_2} = f_{\T_2} + b_{1} f_{\T_1} $. For the right hand side
$ (a_{2   } x_{\T_0} +  a_{   3} x_{\T_1}) T_{i_2} $ 
we have from {\eqref{action}} an expansion of the form
$ (a_{2   } x_{\T_0} +  a_{   3} x_{\T_1}) T_{i_2}=  \sum_{\s \in \tab(\lambda)} c_{\s} x_{\s}  $
Combining, and using the first step of the algorithm, we get 
\begin{equation}\begin{array}{l}
f_{\T_2} + b_{1} f_{\T_1} = \sum_{\s \in \tab(\lambda)} c_{\s} x_{\s}  \, \, \, \, \,  \Longleftrightarrow \\ [.3cm]
f_{\T_2} = -b_{1} (a_{2} x_{\T_0} +  a_{3} x_{\T_1}   )  + \sum_{\s \in \tab(\lambda)} c_{\s} x_{\s}   =
\sum_{\s \in \tab(\lambda)} d_{\s} x_{\s}   
\end{array}
\end{equation}
for some $ d_{\s} \in \cal K $. This is the second step of the algorithm. The following steps of the algorithm
are essentially identical to the second step: as mentioned
above $ f_{\T_{j-1}} T_{i_j} $ is always calculated using the third case of
{\eqref{Seminormal_form}}. The $ k$'th step of the algorithm gives the promised expansion of $ f_\T$.

\medskip
{\bf We call this algorithm for calculating $ f_\T $ 'repeated use of YSF'.}

\medskip
\noindent
{\bf Example}. 
Let $ \lambda := (2,1,1) $ and let $ \T := \T_{\lambda} $, that is

\begin{equation}\label{tableau1}
\T =   \raisebox{-.5\height}{\includegraphics[scale=0.7]{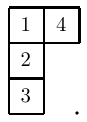}}\, \, \, 
\end{equation}
{\color{black}We have $ d(\T) = s_2 s_3 $ and}
the associated series of $\lambda$-tableaux $ \T_0, \T_1, \T_2 $ is then as follows

\begin{equation}\label{tableaux2and3A}
\T_0 = \T^{\lambda}=  \raisebox{-.5\height}{\includegraphics[scale=0.7]{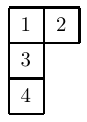}},\, \, \, \, \, \, 
\T_1 =   \raisebox{-.5\height}{\includegraphics[scale=0.7]{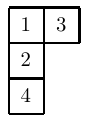}},\, \, \, \, \, \, 
\T = \T_2 =   \raisebox{-.5\height}{\includegraphics[scale=0.7]{tableau1}}\, \, \, 
\end{equation}
Step 1 of the algorithm gives 
\begin{equation}
f_{\T_1} = x_{\T_1} +\dfrac{1}{[2]_q} x_{\T_0}. 
\end{equation}  
Step 2 of the algorithm, using {\eqref{action}}, then gives 
\begin{equation}{\label{regne}}
  \begin{array}{l}
  f_{\T_2} = \dfrac{1}{[3]_q}\left(x_{\T_1} +\dfrac{1}{[2]_q} x_{\T_0}\right) + \left(x_{\T_1} +\dfrac{1}{[2]_q} x_{\T_0}\right) T_{3} \Longleftrightarrow \\[.3cm]
  f_{\T} =   f_{\T_2} =   \dfrac{1}{[3]_q}\left(x_{\T_1} +\dfrac{1}{[2]_q} x_{\T_0}\right) +
  x_{\T_2} +\dfrac{q}{[2]_q} x_{\T_3}  + \dfrac{q-1}{[2]_q} x_{\T_0} 
  \end{array}
  \end{equation}  
where $ \T_3 $ is the tableau
\begin{equation}\label{tableaux2and3}
\T_3 :=  \raisebox{-.5\height}{\includegraphics[scale=0.7]{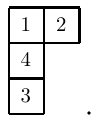}}\, \, \, \, \, \, 
\end{equation}

\medskip
Thus even for this small example the algorithm produces nonstandard tableaux
that must be straightened with the Garnir relations, {\eqref{Garnir}}. This straightening procedure 
is in general a
complicated combinatorial procedure that often has to be repeated many times
until arriving at the desired linear combination
of standard tableaux. In the following sections we shall however see that in the context
of the above algorithm,
the straightening 
procedure can be controlled, at least in certain nontrivial cases.

\section{A generalization of the James-Murphy tableaux }\label{A generalization of the James-Murphy tableaux}
In this section we explain the kind of tableaux in which we shall be particularly interested.

\medskip
For $ 1 \le i \leq j \le n  $ we introduce $ \sigma_{ij} \in \Sy_n $ 
as follows
\begin{equation}
  \sigma_{ij} := \left\{ \begin{array}{ll} s_i s_{i+1}  \cdots  s_{j-1}  & \mbox{ if } i < j \\ 1 &
     \mbox{ if } i = j. \end{array} \right. 
\end{equation}
In other words, $ \sigma_{ij} $ is a one-cycle permutation.
We extend this notation to the Hecke algebra via $ T_{ij} := T_{\sigma_{ij}}$.

\medskip
Let us now 
fix $ \lambda \in \Par_n $ and $ 1 \leq a \leq n $. Suppose
that the $ a$-node of $ \T^{\lambda} $, that is $ (\T^{\lambda})^{-1}(a) $, is a {\it removable}
node of $ {\cal Y}(\lambda) $, meaning that when it is removed from 
$ {\cal Y}(\lambda) $ we still get 
the Young diagram of a partition. 
For any $ a \le b \le n $, we define the tableau $ \T_b $, as follows
\begin{equation}{\label{standardformula}}
\T_b  := \T^{\lambda} \sigma_{ab}. 
\end{equation}  
Then $ \T_b $ is always a standard tableau. 
For example, for $\lambda=(6,6,5,3,3,3) \in \Par_{26}$ and $ a= 12$, we have
$ \T^{\lambda} $, $ \T_{19} $ and $ \T_{25} $ as follows

\begin{equation}\label{moreexamplesTableaux_b}
\T^{\lambda} = \raisebox{-.5\height}{\includegraphics[scale=0.7]{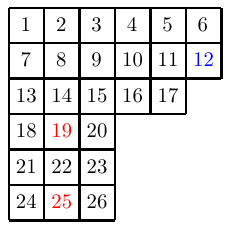}}, \, \, \, \, \, \, \, \, \, \, \, \,
  \T_{19} =  \raisebox{-.5\height}{\includegraphics[scale=0.7]{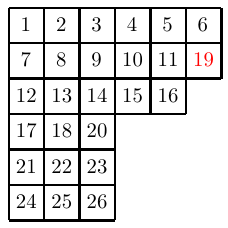}}, \, \, \, \, \, \, \, \, \, \, \, \,
  \T_{25} =  \raisebox{-.5\height}{\includegraphics[scale=0.7]{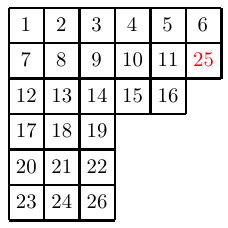}}.
\end{equation}

We shall be especially interested in the cases where 
the $ b $-node of $ \T^{\lambda}$
belongs
to the {\it right border} of $ {\cal Y}(\lambda) $, that is $ b $ and $ b+1 $ lie in different rows of 
$ \T^{\lambda} $ or $ b = n $. For example, keeping $ \lambda = (6,6,5,3,3,3) $, 
the values of $ b $ for which the $b$-node of $ \T^{\lambda}$ belongs to the right border of
$ {\cal Y}(\lambda) $
are $ b = 17, 20, 23, 26$. Below we depict $ \T^{\lambda} $ with
$ a $ in blue, and $ b \neq a $ in red
\begin{equation}\label{examplesTableaux_bA}
\T^{\lambda} = \T_{12}= \raisebox{-.5\height}{\includegraphics[scale=0.7]{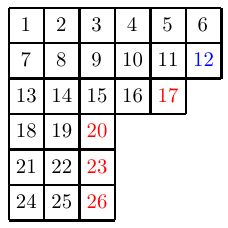}}.\, \, \, \, \, \, 
\end{equation}

The corresponding tableaux $ \T_b $ are as follows
\begin{equation}\label{examplesTableaux_b}
  \T_{17} =  \raisebox{-.5\height}{\includegraphics[scale=0.7]{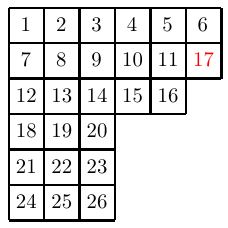}}, \, \, \, \, \, \,
  \T_{20} =  \raisebox{-.5\height}{\includegraphics[scale=0.7]{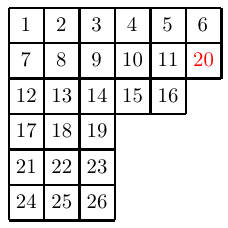}},\, \, \, \, \, \,
  \T_{23} =  \raisebox{-.5\height}{\includegraphics[scale=0.7]{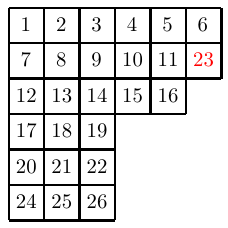}},\, \, \, \, \, \,
      \T_{26} =  \raisebox{-.5\height}{\includegraphics[scale=0.7]{tableau9}}\, \, \, \, \, \,
\end{equation}

\medskip
In the particular case where $ b = n $, the tableau $ \T_b = \T_n $ has the following property
\begin{equation}
 \T_n|_{ 1,2, \ldots, n-1 }  = \T^{\mu } 
\end{equation}  
where $ \mu := shape( \T_n|_{ 1,2, \ldots, n-1 })$.
This kind of tableaux plays an important role in 
James and Murphy's paper \cite{JM}, where the determinant of 
$ \langle \cdot, \cdot \rangle_{\lambda} $ is calculated, 
and for this reason we call the $ \lambda$-tableaux of the form $ \T_b$ {\it generalized James-Murphy tableaux}.
In \cite{JM}, the authors 
prove a formula
for $ \langle f_{\T_n},  f_{\T_n}  \rangle_{\lambda} $ that gives rise to a {\it branching rule} for 
$ \langle \cdot, \cdot \rangle_{\lambda} $ over $ \cal K $. From this branching rule 
they calculate the determinant itself by induction.

\section{YSF along one row}\label{YSF along one row}
In this section we shall see that when we 
work along a single row of $ \T^{\lambda}$, in a sense that we shall shortly explain, 
there is a simple relationship between the corresponding seminormal elements, 
involving only
one denominator.

\medskip
For $ b \ge a $, we consider the tableau $ \T_b $ as in the 
previous section and write for simplicity 
\begin{equation} x_b := x_{\T_b}, \, \,  \,\, \, \, \,  \, \,  \,\, \, \, \,  f_b := f_{ \T_b},
  \, \,  \,\, \, \, \,  \, \,  \,\, \, \, \,
  c_b(i) := c_{ \T_b}(i).
  \end{equation}
We also define $ r_b $ via
\begin{equation}\label{wedefiner_i}
  r_b :=  c_{\lambda}(a)- c_{\lambda}(b). 
\end{equation}
Let $ \{ b_0, b_1, \ldots, b_m\}$ be the 
values of $ b $ corresponding to the right border of $ \T^{\lambda} $, ordered increasingly with $ b_0 = a $,
as in \eqref{examplesTableaux_bA} and \eqref{examplesTableaux_b}.
We are interested in the relation between $ f_{b_i}$ and $ f_{b_{i+1}}$.
We have the following Lemma, that shall be used throughout.

\begin{lemma}{\label{first_lemma}} 
With the above notation the following formula holds 
\begin{equation}\label{thesumofhecke}
  f_{{b_{i+1}}} = f_{{b_i}} \left(
  T_{b_i, b_{i+1}} + \dfrac{1}{\, [r_{b_{i+1}}]_q}\,
  (1+ T_{b_i, b_{i}+1}     + T_{b_i, b_{i}+2}  +  \ldots + 
  T_{b_i, b_{i+1}-1}
  )
  \right).
\end{equation}
(Note that the occurring $ T_{b_i, \beta} $'s have the second index $\beta$ running over   
the row of $ \T_{b_i} $ that contains $ b_{i+1}$). 
\end{lemma}
\begin{dem}
Let $ c_i:= c_{{b_i}}(b_i)- c_{{b_i}}(b_i+1)    $ be the radial distance in $ \T_{b_i} $ from the $ b_i $-node  
to the $ (b_i+1) $-node. Note that $ c_i = c_{\lambda}(a) - c_{\lambda}(b_i+1) $, see for example
the tableaux 
in  \eqref{examplesTableaux_b} of the previous section.
The $ (b_i+1) $-node is situated in a row below the $ b_i $-node in $  \T_{b_i} $
and so when applying YSF we are in the third case of {\eqref{Seminormal_form}}, that is 
\begin{equation}{\label{found}}
f_{{b_i}}T_{b_i}  = f_{{b_i+1}} - \dfrac{1}{\, [c_i]_q } f_{{b_i}} \Longleftrightarrow 
f_{{b_i+1}} =   f_{{b_i}}     \left( T_{b_i}  + \dfrac{1}{\, [c_i]_q } \right).
\end{equation}
We now continue with $ f_{ ({b_i+1}) }T_{b_i+1} $. 
The radial distance in $ t_{ {b_i+1}  } $
from the $ b_i+1 $-node to the $ b_i +2 $-node is $ c_i -1   $
and so we get, using YSF as before, that  
\begin{equation}\label{fwithdeno}
f_{ {{b}_i+1} }T_{{b}_i+1}  = f_{{{b}_i+2}} - \dfrac{1}{\, [c_i-1 ]_q } f_{{{b}_i+1}}  \Longleftrightarrow  
f_{{{b}_i+2}} = f_{ {{b}_i+1} }T_{{b}_i+1}  + \dfrac{1}{\, [c_i-1 ]_q } f_{{{b}_i+1}}.
\end{equation}
We combine this with the expression 
for $ f_{{b_i+1}} $ found in {\eqref{found}} and get 
\begin{equation}
\begin{array}{l}
f_{{{b}_i+2}} = 	
f_{{{b}_i}} 
\left( T_{{b}_i}  + \dfrac{1}{\, [c_i]_q } \right) T_{{b}_i+1}
+
f_{{{b}_i}} \left( T_{{b}_i}  + \dfrac{1}{\, [c_i]_q } \right) \dfrac{1}{\, [c_i-1 ]_q } = \\[.3cm]
f_{{{b}_i}} \left( T_{{b}_i,{b}_i+2} + \dfrac{q}{\, [c_i]_q } +
\dfrac{1}{\, [c_i-1 ]_q } \, T_{{b}_i}
+ \dfrac{1}{\, [c_i ]_q } \dfrac{1}{\, [c_i-1 ]_q } \right) = \\[.3cm]
f_{{{b}_i}}    \left( T_{{b}_i,{b}_i+2}  +
\dfrac{1}{\, [c_i-1 ]_q } \, T_{{b}_i}
+  \dfrac{1}{\, [c_i-1 ]_q } \right)  = 
f_{{{b}_i}}     \left( T_{{b}_i,{b}_i+2}  +
\dfrac{1}{\, [c_i-1 ]_q } \, ( T_{{b}_i} + 1 )  \right)  
\end{array}
\end{equation}
where for the second equality we used $ f_{{b_i}} T_{b_i+1} = qf_{{b_i}} $, which is the first case of
{\eqref{Seminormal_form}}, and for the third equality 
the following Gaussian integer identity
\begin{equation}\label{gaussian}
  q [k-1]_q +1 = [k]_q
\end{equation}
with $ k=c_i$. For the denominator $ [c_i -1]_q $ appearing in the final expression in
\eqref{fwithdeno} we have that 
\begin{equation}
[c_i -1]_q  = [c_{\lambda}(a) - c_{\lambda}(b_i+2)]_q
\end{equation}

This argument is repeated for $ \T_{b_i+3}, \T_{b_i+4},\ldots $ and so on, until
$ \T_{b_{i+1}}$. When we arrive at $ \T_{b_{i+1}}$ the expression for
$f_{b_{i+1}} $ will be as in \eqref{thesumofhecke}
and the 
denominator involved 
will be
\begin{equation}  
  [c_{\lambda}(a) - c_{\lambda}(b_{i+1})]_q = [r_{b_{i+1}}]_q
  \end{equation}
which proves the Lemma.
\end{dem}

\medskip
In the setting of \eqref{examplesTableaux_bA} and \eqref{examplesTableaux_b} the Lemma
gives for example
\begin{equation}
\begin{array}{l}
  f_{ \T_{17}} = f_{ \T_{12}} \left(  T_{12, 17} + \dfrac{1}{[2]_q}(
1+ T_{12, 13} + T_{12, 14} + T_{12, 15} +
  T_{12, 16}
  ) \right) \mbox{ and }   \\[.3cm]
f_{ \T_{20}} = f_{ \T_{17}} \left(  T_{17, 20} + \dfrac{1}{[5]_q}(1 + T_{17, 18} + T_{17, 19}
) \right).
\end{array}
\end{equation}

\medskip
One should observe that it is not possible to have $ r_1 = 1 $ in the Lemma
since that would imply that the $ a $-node of $ \T^{\lambda}$ is not removable,
contrary to the hypothesis of the Lemma. Even so we could still define 
$ \T_{b_1}$ using the same
formula {\eqref{standardformula}}, 
but it would be a nonstandard tableau. 
Note that the sum 	
of Hecke algebra elements in 
\eqref{thesumofhecke}
in this 'limiting case' is exactly the same as the sum of Hecke algebra elements of 
a Garnir relation. This coincidence lies at the heart of 
the results of the following sections.

\medskip
We should remark that Theorem 5.5 of \cite{Ram} offers an alternative approach 
to Lemma {\ref{first_lemma}}. Ram's Theorem 
relates the Garnir relations to Young's seminormal form, in fact 
in the general setting 
of calibrated representations for affine Hecke algebras.
We should however also point out that 
the arguments of the forthcoming sections only rely on one specific Garnir relation,
namely the one mentioned above, 
and hence the results of \cite{Ram} do not offer 
an alternative approach to the rest of our article.

\section{Fat hook partitions}\label{Fat hook partitions}
We assume in this section that 
$ \lambda $ is 
{\it a fat hook partition}, by which we mean that it has the form 
$$ \lambda = ( \lambda_1^{k_1}, \lambda_2^{k_2} ) := 
( \, \stackrel{k_1}{\overbrace{\lambda_1, \ldots  ,\lambda_1}}, \stackrel{k_2}{\overbrace{\lambda_2,\ldots,\lambda_2}} \, ), $$
where $ k_1, k_2 \geq 1$.
Then $ \lambda \in \Par_n $ with $ n = k_1 \lambda_1 + k_2 \lambda_2 $ and 
$ \lambda $ has exactly two removable nodes. We focus on the rightmost of these, the one with coordinates
$ [k_1, \lambda_1 ] $, and 
set $ a := t^{\lambda}[k_1, \lambda_1 ] $. 
As an example we take $ \lambda = (6^2,4^3) $ where $ a = 12 $, that is 
\begin{equation}\label{fattableaux}
  \T^{\lambda} =  \raisebox{-.5\height}{\includegraphics[scale=0.7]{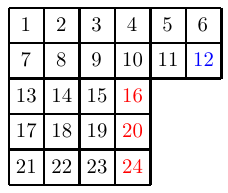}} \, \, \, \, \, \, \, \, \, \, \, \,
\end{equation}

\medskip
For $ i =1,2,\ldots, k_2 $, we associate with 
the $ i$'th row of $ \T^{\lambda} $, counted from
below $ [k_1, \lambda_1] $, an element $ R_i $ of $ {\cal H}^{\cal A}_n(q) $ 
as follows. Let $ \{b_0, b_1, b_2, \ldots, b_{k_2} \} $ denote the 
right border elements of $ \T^{\lambda}$ as before. Then 
the $ i$'th row of $ \T^{\lambda} $, counted from $ [k_1, \lambda_1] $, contains the numbers 
$ \{ b_{(i-1)}+1, b_{(i-1)}+2 , \ldots , b_{i}-1, b_i \} $. Partially inspired
by Lemma \ref{first_lemma}, 
we then
define 
$R_i \in {\cal H}^{\cal A}_n(q) $ via
\begin{equation}{\label{R_i} }
R_i := 1    +  T_{b_{(i-1)}, b_{(i-1)}+1   }  + T_{b_{(i-1)},b_{(i-1)}+2 } + \ldots      + T_{b_{(i-1)},b_i -1} 
\end{equation}
(thus $ b_i $ is skipped).
For example, with $ \lambda $ as in \eqref{fattableaux} we have
\begin{equation}
R_1=1 + T_{12,13} +T_{12,14} +T_{12,15}   \mbox{ and } R_2:= 1+ T_{16,17} + T_{16,18} + T_{16,19}.
\end{equation}  
We next define $ F_{1} := x_{\lambda}  \,R_{1} = x_{b_0}  \,R_{1}   $
and recursively for $ i= 2 , \ldots ,  k_2 $ 
\begin{equation}{\label{fat_hook_recursion}}
F_{i} := 
 ( x_{b_{i}}  - q F_{i-1}   ) R_{i}.
\end{equation}
For example, with $ \lambda $ as in \eqref{fattableaux} we have 
\begin{equation}
F_3=(x_{b_2}-q(x_{b_1}-qx_{b_0}R_{1})R_2)R_3 =(x_{20}-q(x_{16}-qx_{12}R_{1})R_2)R_3.
\end{equation}  
Let us consider the expansion of $ F_{k_2} $ in terms of $x_\T$'s. 
We first observe the 
following useful reformulation of $ F_{k_2} $
\begin{equation}\label{hahaha}  
F_{k_2} = \sum_{j= 1}^{k_2} (-q)^{k_2-j} x_{b_{j-1}} R_{j} R_{j+1} \cdots R_{k_2}
\end{equation}
that follows directly from the definitions. Note that the $ R_j$'s commute.
Now choosing any summand $ T_{\sigma^j} $
from each $ R_j $, one checks easily from \eqref{hahaha} that 
\begin{equation}\label{onecheckseasily}
  x_{b_{j-1}} T_{\sigma^j}  T_{\sigma^{j+1}} \, \cdots \,  T_{\sigma^{k_2}}  =
  x_{{b_{j-1}} \sigma^{j+1} \sigma^{j+2} \,  \ldots \, \sigma^{k_2}  } 
\end{equation}  
and that this is a standard element. 
From this we deduce that the expansion of $ F_{k_2} $ consists of standard elements.

\medskip
With the notation of \eqref{wedefiner_i} 
we set $ r := r_{n} $. For example, in the above case \eqref{fattableaux} we
have that $ r = 5$. 
Our first Theorem is the following surprising formula for 
the expansion of $ f_{n} $ in terms of standard elements, involving only one denominator.
\begin{Thm}{\label{fat_hook_lemma}} We have $ f_n = e_n + \frac{1}{\, [r]_q} F_{ k_2} $. 
  Moreover, the expansion of $ F_{k_2} $ gives rise to a linear combination of standard $ x_\T \! $'s
  as explained above.
\end{Thm}
\begin{dem}
  The second statement was proved in \eqref{onecheckseasily} so let us concentrate on the
  first statement, that is the formula for $ f_n$.
We prove by induction on $ j $ that 
\begin{equation}{\label{inductive_step}}
 f_{b_j} = x_{b_j} + \frac{1}{\, [{\color{blue} r_{b_j}}]_q } F_{ j}  
\end{equation}
from which the formula follows by setting $ j = k_2  $.
The colour blue is only meant to help visualizing the cancellations that take place. 
The induction basis $ j= 1 $ follows directly from Lemma \ref{first_lemma} and the definitions, 
so let us prove the induction step from $ j-1 $ to $j$.
Thus we assume that 
$ f_{b_{j-1}} = x_{b_{j-1}} + \frac{1}{\, [{\color{blue}  r_{b_{j-1}} -1 }]_q }\,  F_{ j-1}. $
From this we deduce via Lemma \ref{first_lemma} that 
\begin{equation}{\label{0}}
\begin{array}{rl}
f_{b_{j}} = & (x_{b_{j-1}} + \frac{1}{\, [{\color{blue} r_{b_j} -1}]_q } F_{ j-1}) T_{ b_{j-1}, b_j} 
+ \frac{1}{\, [{\color{blue}r_{b_j}} ]_q } (x_{b_{j-1}} + \frac{1}{\, [{\color{blue}r_{b_j} -1}]_q }
F_{ j-1}) R_{ j}  \\[.3cm] = &
x_{b_{j}} +  
\frac{1}{\, [{\color{blue}r_{b_j}} ]_q } x_{b_{j-1}} R_{ j} +
\frac{1}{\, [{\color{blue}r_{b_j} -1}]_q } (F_{ j-1} T_{ b_{j-1}, b_j}
+ \frac{1}{\, [{\color{blue}r_{b_j} }]_q }   F_{ j-1} R_{ j} ) .
\end{array}
\end{equation}
We consider the two terms of the last parenthesis. 
Let $ x_\T$ be a standard element occurring in the expansion of 
$ F_{ j-1} $, in the sense explained in \eqref{onecheckseasily}. Then the action of $ R_{j} $ on  
$ x_\T$ only involves the 'easy' second case of {\eqref{action}}. To be more precise,
let $ \sigma^1, \sigma^2, \ldots, \sigma^{\lambda_2}$ be the permutations
involved in the definition of $ R_{j} $, corresponding to the numbers
of the $ j$'th row of $ \T^{\lambda}$ counted from $ [k_1, \lambda_1] $, see
\eqref{R_i}. Then we have that
\begin{equation}\label{wekeep}
  x_\T R_{j} = x_{ \T \sigma^1} + x_{ \T \sigma^2} + \ldots + x_{ \T \sigma^{\lambda_2}}. 
\end{equation}
We 
keep $ x_{\T} $, but now 
focus on the $ x_\T T_{ b_{j-1}, b_j} $ term of the last parenthesis of {\eqref{0}}.
Recall that by definition 
$ T_{ b_{(j-1)}, b_j} = T_{b_{(j-1)} } T_{ b_{(j-1)}+1 }  \cdots  T_{ b_{j}-1} $. Using this we 
get once again via {\eqref{action}} that  
\begin{equation}{\label{2}}
x_\T  T_{ b_{(j-1)}, b_j} =  x_{ \T \sigma_{b_{(j-1)}, b_j}}
\end{equation}
but this $ x_{\T \sigma_{ (j-1), b_j }}  $ is not a standard element.
To illustrate,  we take $ \lambda $ as in \eqref{fattableaux}, $ j= 2$, and 
$ x_\T$ occurring in the expansion of $ F_{ j-1}= F_{1} $ with $ \T $ and 
${ \T \sigma_{b_{(j-1)}, b_j}} =  \T \sigma_{16, 20}  $ as follows
\begin{equation}\label{fattableauxA}
  \T =  \raisebox{-.5\height}{\includegraphics[scale=0.7]{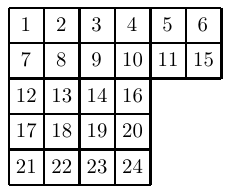}}, \, \, \, \, \, \, \, \, \, \, \, \,
  \T  \sigma_{16, 20} =  \raisebox{-.5\height}{\includegraphics[scale=0.7]{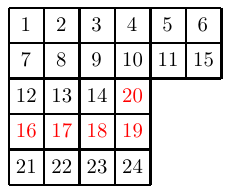}}. 
\end{equation}
The numbers where $ \T $ and $ \T \sigma_{16, 20} $ differ have been coloured red, 
they form a Garnir belt.

\medskip
The last comment holds in general. 
Let therefore $ \g_{k_1+j {\color{black}{-1}}, \lambda_2} $ be the Garnir $ \lambda$-tableau, as introduced in the section above
\eqref{garnir}. Then we have that $ d(\g_{k_1+j{\color{black}{-1}}, \lambda_2}) = \sigma_{b_{(j-1)}, b_j}  $
which is a $ \Sy_n $-element commuting with $ d(\T)$
and so the Garnir relation {\eqref{Garnir}}
gives that 
\begin{equation}\label{garnirgives}
  x_{\T \sigma_{ (j-1), b_j }} + x_{ \T \sigma^1} + x_{ \T \sigma^2} + \ldots + x_{ \T \sigma^{\lambda_2}} =0
\end{equation}
where the $ \sigma^j $'s are as in \eqref{wekeep}. 
Combining \eqref{garnirgives}, \eqref{0} and \eqref{2}, and
using once again the Gaussian identity \eqref{gaussian}, we arrive at 
\begin{equation}{\label{4}}
\begin{array}{l}
f_{b_{j}}  = 
x_{b_{j}} +  
\frac{1}{\, [{\color{blue}r_{b_j}} ]_q } x_{b_{j-1}} R_{ j} -
\frac{q}{\, [{\color{blue}r_{b_j}} ]_q } 
  F_{ j-1} R_{ j} = \\ [.3cm]
x_{b_{j}} +  
\frac{1}{\, [{\color{blue}r_{b_j}} ]_q } ( x_{b_{j-1}}  -
q F_{ j-1}) R_{ j}  =  x_{b_{j}} + \frac{1}{\, [{\color{blue}r_{b_j}} ]_q } F_{ j}. 
\end{array}
\end{equation}
This completes the inductive step of the proof of the Lemma.
\end{dem}

\medskip \medskip
Let us illustrate the formula on the partition $ \lambda = (3, 2^2 ) $ of $7 $. In that
case we have $ r = 3 $ and the formula for $ f_7 $ becomes 
\begin{equation}\label{fattableauxB}
 \raisebox{-.4\height}{\includegraphics[scale=0.7]{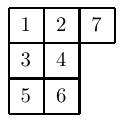}} + \dfrac{1}{[3]_q}
  \left( \raisebox{-.4\height}{\includegraphics[scale=0.7]{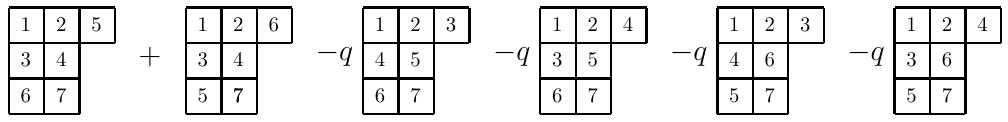}} \right)
\end{equation}
where we identify $ \T $ and $ x_\T $. 

\medskip
Let us now illustrate the Theorem on the bigger example
$ \T_{24}$ for 
$ \lambda = (6^2, 4^3)$, that is
\begin{equation}
  \T_{24} =  \raisebox{-.5\height}{\includegraphics[scale=0.7]{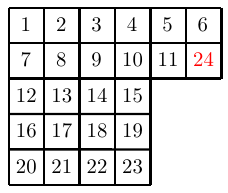}}
\end{equation}
which will also give an indication on how to work with the
Theorem in general, via formula \eqref{hahaha} for $F_{k_2}$.
In this case $ k_2 = 3 $ and so there will be three summands 
$  x_{b_{j-1}} R_{j} R_{j+1} \cdots R_{k_2} $ in \eqref{hahaha} with the following 'leading terms' $  x_{b_{j-1}} $
\begin{equation}\label{leading-terms}
  x_{b_0} =  \raisebox{-.5\height}{\includegraphics[scale=0.7]{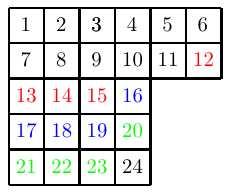}}, \, \, \, \, \, \, 
  x_{b_1} =  \raisebox{-.5\height}{\includegraphics[scale=0.7]{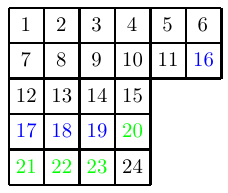}}, \, \, \, \, \, \, 
  x_{b_2} =  \raisebox{-.5\height}{\includegraphics[scale=0.7]{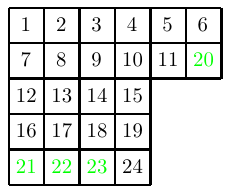}}. 
 \end{equation}
By definition $   x_{b_0} R_1 R_2 R_3 $ is the sum of all the products of
the 'monochromatic' cycles
of  $ x_{b_0} $ where each cycle starts in $ 12, 16 $ or $ 20$, in other words products of a 'red'
cycle $ 1, (12, 13) $,
$(12,13,14) $, 
$ (12,13,14,15) $, a 'blue' cycle $ 1, (16, 17), (16, 17,18), (16, 17,18,19) $ and a 'green' cycle
$ 1$, $ (20,21)$, $ (20,21,22)$ and $ (20,21,22,23) $.
{\color{black}{Note that this can also be described as the sum of all standard elements
obtained from $ x_{b_0} $ by shuffling monochromatic numbers}}. 
Thus $ x_{b_0} R_1 R_2 R_3 $ is the sum of $ 4 \times 3 \times 3 =36  $ 
standard elements, that all enter in $ F_{k_2} $ with coefficient $ q^2 $, and similarly 
$ x_{b_2} R_ 2 R_3 $ is {\color{black}the} sum of $ 12 $ standard elements that all enter in $ F_{k_2} $
with coefficient $ {\color{black}{-}} q $ whereas $ x_{b_3} R_3 $ is {\color{black}the} sum of $ 4 $
standard elements that all enter 
in $ F_{k_2} $ with coefficient $ 1 $.

\medskip
One also checks that the standard elements in $   x_{b_0} R_1 R_2 R_3 $, $   x_{b_1}  R_2 R_3 $ and
$   x_{b_2}  R_3 $ are all different. For example, if $ x_\T$
%{\color{green}{\sout{were}}}
{\color{black}{is}} a standard element appearing in both 
$   x_{b_0} R_1 R_2 R_3 $ and $   x_{b_1}  R_2 R_3 $, then
looking at the distribution of blue numbers in $ x_{b_0}$ and $ x_{b_1}$ one sees that $16$ would have to be in
the fourth row of $ \T $ and therefore either $ 17, 18 $ or $ 19 $
would have to be in the second row of $ \T$, which is impossible.

\medskip
\noindent
%{\color{green}{\sout{
%{\bf Remark}
%Calculations with small partitions indicate that the inverse
%expansion of $ x_\T$ in terms of $ f_\T $ 
%does not permit a 'nice' description as the one of 
%the Theorem.}}}

\section{Expansion of $f_n$ for general partitions}\label{relevantsection}
Our next aim is to show that the results from the previous section can
be extended to arbitrary partitions $ \lambda$. In this section we determine the expansion
of the seminormal element $ f_n$ for the generalized James-Murphy tableaux in terms of standard elements. 
Although this extension does not require substantially new ideas compared to the previous section,
the notational technicalities are 
more involved and so we prefer to treat this extension separately.

\medskip
Let us set up the relevant notation. 
Let $ \lambda $ be a partition of $ n $. 
We fix a removable node $ [\alpha_0, \beta_0]  $ for $ \lambda $ and let 
the removable nodes below
%{\color{green}{\sout{$ [\alpha_0, \alpha_0]  $}}}
{\color{black}{$ [\alpha_0, \beta_0]  $}}
be $ [\alpha_j, \beta_j],  j = 1, 2 , \ldots , N $, 
from top to bottom. In \eqref{fattableauxC} we give the example $ \lambda = (6^2, 4^3, 3^2, 1) $, 
where we choose $ [\alpha_0, \beta_0] := [ 2,6] $ and where we
indicate with arrows the $  [\alpha_j, \beta_j] $'s. 
We set $ d_j := \T^{\lambda}[\alpha_j, \beta_j] $ and so have in \eqref{fattableauxC} that
$ d_0 = 12, d_1 = 24, d_2 = 30, d_3= 31$. 
\begin{equation}\label{fattableauxC}
 \raisebox{-.4\height}{\includegraphics[scale=0.7]{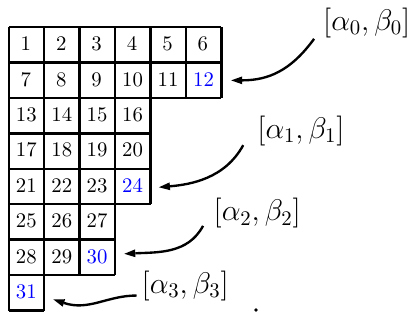}}
\end{equation}

With the notation from the previous sections we have 
$ \T_n =   \T^{\lambda}{ \sigma_{d_0,d_N} }$ and our aim is to 
determine the expansion of 
$ f_n = f_{\T_n}$ in terms of standard elements $x_\T$. Here is
$ \T_n $ for the same $ \lambda $ as in \eqref{fattableauxC}

\begin{equation}\label{exJMtableau}
\T_n:=   \raisebox{-.5\height}{\includegraphics[scale=0.7]{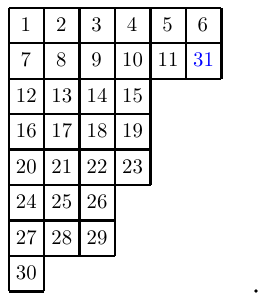}} \, \, \, \, \, \, \, \, \, \, 
\end{equation}

\medskip
Now   
$ \lambda $ determines a series
of subpartitions $ hook_j(\lambda), j = 1, 2, \ldots, N $ of $ \lambda $ via 
\begin{equation}
 {\cal Y}(hook_j(\lambda)) := shape(\T|_{ 1,2, \ldots, d_j }) \setminus shape(\T|_{ 1,2, \ldots, d_{j-2} })
\end{equation}  
where $ d_{\color{black}{-1}} := 0 $ and $ shape(\T|_{ 1,2, \ldots, 0 }) := \emptyset$.
We then define the $hook_j(\lambda)$-tableau
$ \T^{ hook_j(\lambda)} $ to be the {\it restriction of $ \T^{\lambda} $} to $hook_j(\lambda)$.
Thus the numbers appearing in $ \T^{ hook_j(\lambda)} $ are $ \{ d_{j-2}+1,   d_{j-2}+2, \ldots, d_j \} $
and so, strictly speaking, $ \T^{ hook_j(\lambda)} $ is not a tableau according to the definition in section 2,
but still we shall refer to it as a
$ hook_j(\lambda)$-tableau. Below we give the $ {\color{black}{\T}}^{hook_j(\lambda)}$'s for $ \lambda $
as in \eqref{fattableauxC}.

\begin{equation}\label{morefattableauxD}
\T^{hook_1(\lambda)} =   \raisebox{-.4\height}{\includegraphics[scale=0.7]{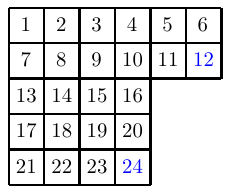}}, \, \, \, \, \, \, \, \, \, \, 
\T^{hook_2(\lambda)} =
\raisebox{-.4\height}{\includegraphics[scale=0.7]{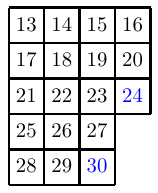}}\, , \, \, \, \, \, \, \, \, \, \, 
\T^{hook_3(\lambda)} =     \raisebox{-.4\height}{\includegraphics[scale=0.7]{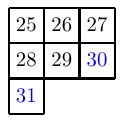}}.
\end{equation}

Note that in general $ hook_j(\lambda) $ is a fat hook partition, except possibly for $ hook_1(\lambda)$.

\medskip
Set $  b_0^j := d_{j-1}$ and 
let $ \{ b_0^j, b_1^j,  \ldots, b_k^j \}$
be the right border of
$ \T^{ hook_j(\lambda)} $ below
the $  b_0^j $-node.
Generalizing \eqref{R_i} we set 
\begin{equation}
R_i^j  :=1  + T_{b_{(i-1)}^j, b_{(i-1)}^j+1}    + T_{b_{(i-1)}^j,b_{(i-1)}^j+2 }+   \ldots +T_{b_{(i-1)}^j,b_i^j -1}.    
\end{equation}
We next define $ {\cal F }_k^j  \in {\cal H}^{\cal A}_n(q)$, generalizing
the element $ F_{k_2}  $ of 
{\eqref{fat_hook_recursion}}.
Set first 
$ {\cal F }^j_1 := R_1^j$ and then recursively
\begin{equation}
{\cal F }^j_{i} := ( T_{b_{0}^j, b_{i}^j}-     q {\cal F }^j_{i-1} ) R_i^j. 
\end{equation}
Finally set
\begin{equation}
{\cal F }^j := {\cal F }^j_k. 
\end{equation}
Generalizing \eqref{hahaha}, we have the following reformulation of 
$ {\cal F }^j  $ 
\begin{equation}\label{using}
  {\cal F }^j =  \sum_{i= 1}^{k} (-q)^{k-i} T_{b_0^j, b_{i-1}^j} R_{i}^j R_{i+1}^j \cdots
  R_{k}^j. 
\end{equation}
Note that there is difference between $ {\cal F }^j  $ and $ F_k $: the former belongs 
to $ {\cal H}^{\cal A}_n(q) $ whereas the latter belongs to $ S_q^{\cal A}(\lambda)$.

\medskip
\noindent
To illustrate these definitions we 
write down the $ R_i^j$'s and $ {\cal F}^j$'s, for $ \lambda $ as in \eqref{fattableauxC}.
Using \eqref{using} we get
\begin{equation}\label{toillustratethesedef}
\begin{array}{l}
  R_1^1= 1+T_{12,13}+  T_{12,14}+T_{12,15},  \, \, 
  R_2^1= 1 +T_{16,17}+  T_{16,18}+T_{16,19} \, \,   \\ [.3cm]  
    R_3^1= 1 +T_{20,21}+  T_{20,22}+T_{16,23}     \, \,   \\ [.3cm]  
    {\cal F }^1 = q^2     R_1^1    R_2^1 
    R_3^1 - q T_{12,16} \, R_2^1   R_3^1 + T_{12,20}    R_3^1     \\ [.3cm]
    R_1^2= 1+T_{24,25}+  T_{24,26},  
     R_2^2= 1+T_{27,28}+  T_{27,29}      \\ [.3cm]
    {\cal F }^2 =  {\color{black}{-}}q  R_1^2    R_2^2    +  T_{24,27} \, R_2^2     \\ [.3cm]
    R_1^3= 1   \\ [.3cm]
        {\cal F }^3 = 1. 
\end{array}
\end{equation}

\medskip
We have that 
$ r_{d_i} := c_{\lambda}(d_0) - c_{\lambda}(d_i) $, see \eqref{wedefiner_i}. 
We set ${\PP}^0 := 1$ and
recursively for $ j=1,2, \ldots, N$
\begin{equation}\label{mainrecursion}
  {\PP}^i := \PP^{i-1} \left( T_{ d_{i-1}, d_i} + \dfrac{ 1}{[r_{d_i}]_q}{{ \cal F}^i} \right). 
\end{equation}

We are interested in $ x_{\lambda}\, \PP^N \in S_q^{\cal K}(\lambda) $.  
For example, using once again $ \lambda $ as in \eqref{fattableauxC} we get that
\begin{equation}\label{expansionexpansion}
x_{\lambda}  \PP^3 = x_{\lambda}  \left( T_{12,24} + \dfrac{ 1}{[5]_q}     {\cal F }^1  \right)
  \left( T_{24,30} + \dfrac{ 1}{[8]_q }        {\cal F }^2 \right)   \left( T_{30,31} + \dfrac{ 1}{[11]_q}
       {\cal F }^3 \right).
\end{equation}  
where the ${ \cal F}^j$'s are as in \eqref{toillustratethesedef}.

\medskip
The main result of this section, Theorem {\ref{f_n}}, is the identity $ f_n = x_{\lambda}\, \PP^N $. 
We need two auxiliary lemmas. Here is the first one. 
\begin{lemma}{\label{first_auxilliary_lemma}}
For each $j = 1,\ldots, N $ there is a 
$ p_j \in {\cal H}^{\cal K}_n(q)  $ such that 
$  f_{d_{j}} =  f_{d_{j-1}} \,  p_j   $. It satisfies 
$  x_{d_{j-1}} \, p^{j} = x_{d_{j}} + \frac{1}{\, \, [r_{d_j}]_q}  \, x_{d_{j-1}} \,{\cal F}^{j} $.
\end{lemma}
\begin{dem}
  The existence of $ p_j$ follows from repeated applications of Lemma \ref{first_lemma}. 
  It is a product of factors, each one of the form
  \begin{equation}\label{replacingeachfactor}
    T_{b_i, b_{i+1}} + \dfrac{1}{\, [r_{b_{i+1}}]_q}\,  (
1+T_{b_i, b_{i}+1}+T_{b_i, b_{i}+2}     + \ldots +
    T_{b_i, b_{i+1}-1} 
    ). 
\end{equation}  
The formula $  x_{d_{j-1}} \, p_{j} = x_{d_{j}} + \frac{1}{\, \, [r_{d_j}]_q}  \, x_{b_{j-1}} \,{\cal F}^{j} $
follows from arguments identical to those in the proof of 
Theorem \ref{fat_hook_lemma}. Note that these arguments in fact show that
$ f_n = e_{\lambda} \,  p_1 $. The  
cancellations of Theorem 
\ref{fat_hook_lemma} depend only on the Garnir relations and in particular 
they do not require that the nodes above the relevant Garnir belt, which is always of the form as in
\eqref{fattableauxA}, are those of a fat hook partition.
Hence the cancellations will also occur in the present setting. This proves the last
statement of the Lemma.
\end{dem}

\medskip
\medskip
For $ x \in \mathbb Z $ we define $  p^{x}_j  $ the same way as $  p_j  $, but replacing each factor
\eqref{replacingeachfactor} with 
  \begin{equation}\label{replacingeachfactorx}
    T_{b_i, b_{i+1}} + \dfrac{1}{\, [r_{b_{i+1}}+x]_q}\,  (
1+T_{b_i, b_{i}+1}+T_{b_i, b_{i}+2}     + \ldots +
  T_{b_i, b_{i+1}-1} )
  \end{equation}  
and thus we have $ p^{0}_j = p_j $.
We next introduce 
$ f^{x}_j \in  S_q^{\cal K}(\lambda)  $ via
\begin{equation}
f^{x}_j :=    x_{\lambda} \, p^{x}_j. 
\end{equation}
It is 
a generalization of the seminormal basis, in fact
$ f_n = f_{n}^0 $ in the case where $ \lambda $ is a fat hook partition.

\medskip
We need a slightly more general version of this construction.
Let $ \T $ be a $ \lambda$-tableau
that coincides with $ \T^{\lambda}$ in all nodes below the node $ [\alpha_j, \beta_j]$, that is 
\begin{equation}
\T[\alpha,\beta] = \T^{\lambda}[\alpha,\beta] \, \,  \mbox{ if } \alpha > \alpha_j.
\end{equation}
For example, using $ \lambda $ as in \eqref{fattableauxC}, and $ j =2 $, the
following $\lambda$-tableau $ \T $ could be used, since $ {\T}^{\lambda} $ and $ \T $ coincide
below the red line. 

\begin{equation}\label{fattableauxD}
{\T}^{\lambda} := \raisebox{-.5\height}{\includegraphics[scale=0.7]{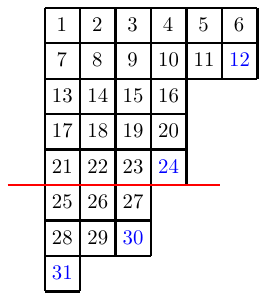}},  \, \, \, \, \, \, \, \, \, \, 
\T:=   \raisebox{-.5\height}{\includegraphics[scale=0.7]{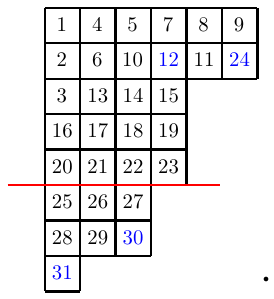}} \, \, \, \, \, \, \, \, \, \, 
\end{equation}

For such a $\lambda$-tableau $ \T $ and $ x \in \mathbb Z $ we define 
\begin{equation}
  f^{x,\T}_j :=  x_\T \, p_j^x.
\end{equation}  
Our second auxiliary Lemma is as follows. 
\begin{lemma}{\label{second_auxilliary_lemma}}
  In the above setup we have
\begin{enumerate} \itemsep 0.2cm
\item   
  $ f^x_{d_{j}} = x_{d_{j}} + \frac{1}{\, \, [r_{d_j}+x ]_q} \,  x_{d_{j-1}} \,{\cal F}^j .$ 
\item  $ f^{x,\T}_{d_{j}} = x_\T T_{\sigma_{d_{j-1},d_{j} }} + 
  \frac{1}{[r_{d_j}+x]_q } \,  x_\T \, {\cal F}^j. $
\end{enumerate}
\end{lemma}
\begin{dem}
Just as in the proof of Lemma {\ref{first_auxilliary_lemma}}, we 
recycle 
the proof of Theorem \ref{fat_hook_lemma}. That proof 
{\color{black}{depends}} on the formula given in Lemma \ref{first_lemma}, and on 
cancellations that arise from the Garnir relations of the form \eqref{fattableauxA}. These cancellations
also take place in the present setting, when we replace the $ \color{blue}{r_{b_j}}$'s of that proof with
$ \color{black}{r_{b_j}+x}$. This proves the Lemma.
\end{dem}
\medskip \medskip

We are now in position to prove the main Theorem of this section, which generalizes  
Theorem \ref{fat_hook_lemma} to arbitrary partitions.
\begin{Thm}{\label{f_n}} 
  Let $   {\PP}^N  \in {\cal H}_q^{\cal K}(\lambda)  $ be
  the element given by the recursion \eqref{mainrecursion}. Then 
  $  x_{\lambda} {\PP}^N  \in S_q^{\cal K}(\lambda) $ satisfies
$ f_n = x_{\lambda} {\PP}^N $. 
Moreover, the $ x_\T$'s arising from the expansion of $ x_{\lambda} {\PP}^N $ 
are all standard elements.
\end{Thm}
\begin{dem}
We proceed by induction on $ N $. The basis case $ N =1 $ is
$ f_{d_1} = x_{d_1} + \frac{1}{[r_{d_1}]_q}  x_{\lambda} \, {\cal F}^1  $, 
and it corresponds to 
Theorem \ref{fat_hook_lemma}.
Let us consider the step from $ N=1 $ to $ N=2$. 
We consider the standard tableaux $ \T $ such that $ x_{\T}$
appear in the expansion of 
$ f_{d_1} = x_{d_1} + \frac{1}{[r_{d_1}]_q}  x_{\lambda} \, {\cal F}^1  $.
In $ \T_{d_1}$, corresponding to the first term
$x_{d_1} $ of this expansion, we have that 
$ d_1 $ is located in
position $ [\alpha_0, \beta_0] $, whereas in all the other $ \T $'s, 
we have
that $ d_1 $ is located in position $[\alpha_1, \beta_1] $.
\medskip

Now 
\begin{equation}
  f_{b_{2}} =  f_{b_{1}} \,  p_2 = x_{d_1}\, p_2 + \frac{1}{[r_{d_1}]_q}  x_{\lambda}  {\cal F}^1 \, p_2.
\end{equation}  
With regards to $  x_{d_1} \, p_2$ we get by 
Lemma \ref{first_auxilliary_lemma} 
that 
\begin{equation}  x_{d_1} p_2 = 
x_{d_{2}} + \frac{1}{[r_{d_2}]_q} \,   x_{d_{1}} { \cal F}^2   = 
x_{d_{1}} \left( T_{d_{1},d_{2}}
+ \frac{1}{[r_{d_2}]_q} \,   { \cal F}^2 \right) 
\end{equation}
and hence we 
must prove that the same formula holds for all $ x_\T $ involved 
in $ \frac{1}{[r_{d_1}]_q}  x_{\lambda}  {\cal F}^1 \, p_2$, 
i.e. that 
\begin{equation}\label{iscorrect}  x_{\T} p_1  = 
  x_{\T} \left( T_{d_{1},d_{2}} + \frac{1}{ [r_{d_2}]_q}   \,   { \cal F}^2 \right)
\end{equation}  
holds for these $ x_\T $. As already mentioned, 
for each such $ x_\T $-term we have that 
$ \T[\alpha_1, \beta_1 ]  = d_1 $. 
We now apply part $b)$ of 
Lemma \ref{second_auxilliary_lemma} 
with $ x = r_2 - r_1 $, and conclude that \eqref{iscorrect} is correct.
The general induction step is treated the same way.
\end{dem}

\medskip
\noindent
%{\color{green}{\sout{
%{\bf Remark}
%We remark that Theorem 3.4 of \cite{FLT2} gives a different approach to 
%the Theorem, which is based on the complete expansion of 
%$ f_n $, (in the symmetric group case). This approach appears to be in the spirit 
%of the expansion explained in \eqref{leading-terms}. 
%It allows the authors of \cite{FLT2} to determine the greatest common divisor of the coefficients
%of the expansion of $ f_n$.}}}

{\color{black}
\medskip
\noindent
We now give a characterization of the standard tableaux that appear
    in the expansion of $ f_n = x_{\lambda} {\PP}^N $ in Theorem \ref{fat_hook_lemma},
    extending equation \eqref{leading-terms} for fat hook partitions and the comments below it.

\medskip
The characterization relies on a set of rules that allow us to produce a set of
standard tableaux $\mathcal{CS}_{\lambda}$, in which all the numbers
are coloured from a large colour set, that includes black. The standard tableaux
that appear in the expansion of $ f_n $ are those tableaux that are obtained
from $\mathcal{CS}_{\lambda}$ by shuffling
in all possible ways the monochromatic numbers, except the black ones, such that the results
are still standard.

\medskip
Let therefore the situation be as above,
with $ \lambda \in \Par_n$ arbitrary, 
having removable nodes $ [\alpha_j, \beta_j],  j = 0, 1, 2 , \ldots , N $, see 
\eqref{fattableauxC}. In the following the words 'before', 'after', 'next' and so on, refer to the
natural total order on 
the nodes of $ \lambda $
given by $ [r_1,c_1] < [r_2,c_2]  $ iff $ \T^{\lambda}[r_1, c_1] <  \T^{\lambda}[r_2, c_2] $; this is
'the row reading order' on the nodes for $ \lambda$. Moreover, the word 'coloured' refers to
a colour different from black.

\medskip
We consider $ \T $ as a  filling of the nodes of $ \lambda $
with the numbers $1, 2, \ldots, n $, in this order.
Then $  \mathcal{CS}_{\lambda} $ is the set of $ \lambda$-tableaux $ \T $ that can be obtained
by applying the following rules. 
\begin{description}
\item[Rule 1] Let $ i := \T^{\lambda}[r,c] $ where $ [r,c] $ is a node strictly 
  before $ [\alpha_0, \beta_0] $. Then $ i $ should also be placed in $[r,c] $ in $ \T $ 
  and should be coloured black.

\item[Rule 2] Let $ i := \T[r,c] $ where 
  neither $ [r,c] $ nor $ [r,c+1] $ belongs to the right border of $ \lambda$. 
  Then $i+1$ should be placed in $ [r,c+1] $ and should be coloured 
with the same colour as $ i $. 

\end{description}
These two rules are visible in \eqref{leading-terms}. 
It follows from them that in $ \T$, 
the node and colour of $ i+1 $ is uniquely determined by the
node and colour of $ i $, unless 
the node of $ i $ is situated 
after $ [\alpha_0, \beta_0] $, and  
belongs either to the right border of $ \lambda $ or is a node one step before the right border of 
$ \lambda$. The remaining rules consider these cases.

\begin{description}
\item[Rule 3] Suppose that $ i = \T[\alpha_0, \beta_0-1 ] $. 
  Then $ i+1 $ should be placed 
  either in $ [\alpha_0, \beta_0 ] $, with a colour different from black, or 
  in $ [\alpha_0, \beta_0 ] $ with colour black. Here an example with $ i =12$. 
\begin{equation}
\raisebox{-.5\height}{\includegraphics[scale=0.7]{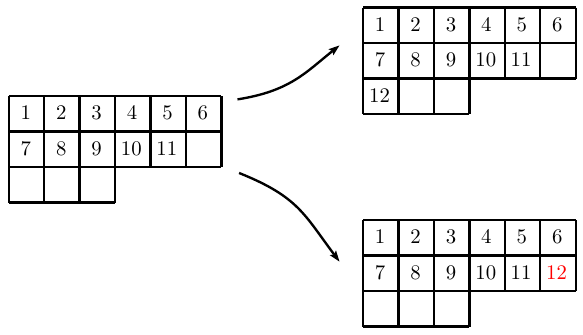}}
\end{equation}

\item[Rule 4] Let $ [r, c] $ be a node one step before the right border of $ \lambda $
  such that $ r > \alpha_0 $  and suppose that $ i:= \T[r,c] $ is black. Then $ i+1 $ should
  be placed in $ [r,c+1] $ and should be black. Here an example with $ i=11$. 
\begin{equation}
\raisebox{-.5\height}{\includegraphics[scale=0.7]{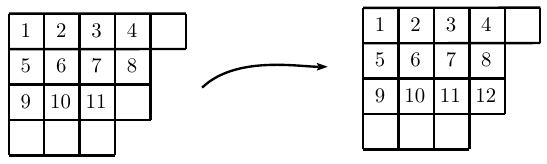}}
\end{equation}

\item[Rule 5] Let $ [r, c] $ be a node one step before the right border of $ \lambda $
  such that $ r \not\in  \{\alpha_0, \ldots, \alpha_N \} $  and suppose that $ i:= \T[r,c] $
  is coloured. 
  Then $ i+1 $ should
  be placed in $ [r,c+1] $ and should have a previously unoccupied colour. Here is an example. 
\begin{equation}
\raisebox{-.5\height}{\includegraphics[scale=0.7]{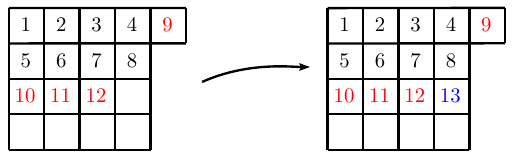}}
\end{equation}

\item[Rule 6] Let $ [r, c] $ be a node one step before the right border of $ \lambda $
  such that $ r \in   \{\alpha_1, \ldots, \alpha_N \} $  and suppose that $ i:= \T[r,c] $ is coloured.
  Then $ i+1 $ should either 
  be placed in $ [r,c+1] $ and have a previously unoccupied colour or should be
  placed in $ [r+1,1] $ with colour black. 
\begin{equation}
\raisebox{-.5\height}{\includegraphics[scale=0.7]{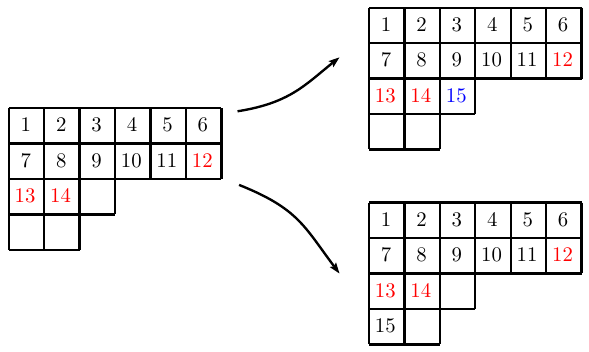}}
\end{equation}

\item[Rule 7]
Suppose that $ i := \T[r, c]  $ is black and that 
$ [r, c] $ is a node belonging to the right border of $ \lambda $
  such that $ r > \alpha_0 $. Then $ i+1 $ should be placed
  either in the first node of the $ r+1$'st row of $ \lambda$, with colour black, or 
  in the unique previously unoccupied node $ [\alpha_i, \beta_i] $ where $ \alpha_i < r$, with a 
  previously unoccupied colour, in particular different from black. Here is an illustration with $ i=12$.
\begin{equation}
\raisebox{-.5\height}{\includegraphics[scale=0.7]{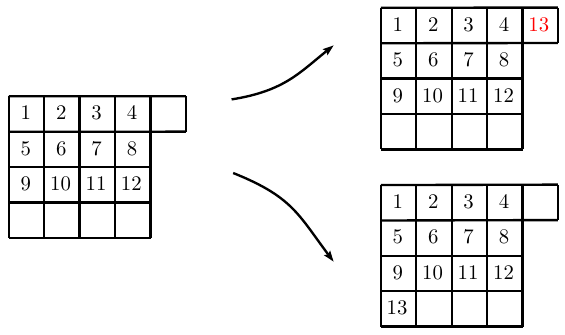}}
\end{equation}

\item[Rule 8] Suppose that $i:= \T[r,c] $ where $ [r,c] $ is a
  node belonging to the right border of $ \lambda $, such
  that $ r \ge \alpha_0 $ and such that $i$ is coloured. Then $ i+1 $ should be placed 
  in the first node of the first unoccupied row of $ \lambda$, with the same colour as $ i$,
  unless that node belongs to the right border of $ \lambda $ in which case it should have a new colour. 
  We illustrate with $ i =13$. 
\begin{equation}
\raisebox{-.5\height}{\includegraphics[scale=0.7]{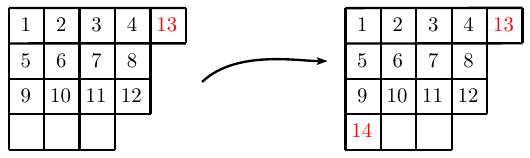}}
\end{equation}

\end{description} 
Observe that applying these rules, filling in a black $i $ in the $r$'th row
where $ r > \alpha_0$ there
will always be a unique unoccupied node before $ i$, whereas filling in a coloured $ i $ in
row $ r $, there will be no such unoccupied node.

\medskip
The rules can be read off directly from \eqref{mainrecursion}. Indeed, the black
numbers of the rules correspond to the factors $ T_{ d_{i-1}, d_i} $ of \eqref{mainrecursion} and
the coloured numbers correspond to the factors ${ \cal F}^i $ of \eqref{mainrecursion}.
The connection with
\eqref{mainrecursion}
also allows us to the determine the coefficient of $ x_\T $ in $ f_n$, since black numbers correspond 
to $ T_{ d_{i-1}, d_i} $, for some $ i $, 
that only contributes with 1 to $f_n $, whereas $ j $ coloured numbers that correspond to 
$ { \cal F}^i $, for some $ i $, contribute
with $  \dfrac{(-q)^{j-1}}{ [r_{d_i}]_q } $ to $ f_n$. Note that for standard tableaux $ \s $ and $ \T$
related via a shuffling of monochromatic numbers, the coefficients of
$ x_\s $ and $ x_\T $ are the same.

\medskip
Setting $ \mu := \lambda \setminus [\alpha_0, \beta_0] $, 
the set $ \mathcal{CS}_{\lambda} $ 
also appears in \cite{FLT2}
where it is denoted the set of {\it colour semistandard tableaux} $ SStd(\mu; 1 ) $, 
although the formulation in \cite{FLT2} is different from ours.

\medskip
Suppose for example that $ \lambda= (4,3,2,2) $ with $ (\alpha_0, \beta_0 ) =(1,4)$. Then
applying the rules we get the following 6 elements of $ \mathcal{CS}_{\lambda}$, with corresponding
coefficients.  
\begin{equation}\label{sixexample}
  \mathfrak{s}_0 = \raisebox{-.6\height}{\includegraphics[scale=0.7]{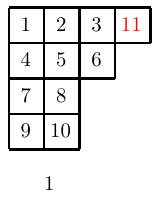}} \, \,  
  \mathfrak{s}_1 = \raisebox{-.6\height}{\includegraphics[scale=0.7]{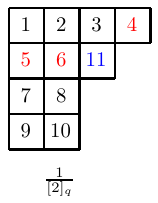}}\, \, 
  \mathfrak{s}_2 = \raisebox{-.6\height}{\includegraphics[scale=0.7]{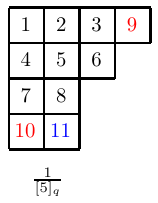}}\, \, 
  \mathfrak{s}_3 = \raisebox{-.6\height}{\includegraphics[scale=0.7]{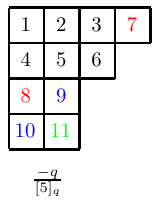}}\, \, 
  \mathfrak{s}_4 = \raisebox{-.6\height}{\includegraphics[scale=0.7]{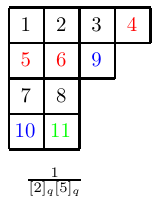}}\, \, 
  \mathfrak{s}_5 = \raisebox{-.6\height}{\includegraphics[scale=0.7]{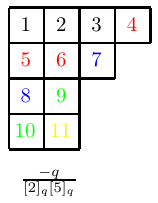}} \, 
\end{equation}
This is the 'running example' of \cite{FLT2} 
and one sees that the tableaux $   \mathfrak{s}_0, \ldots,   \mathfrak{s}_5 $ in
\eqref{sixexample} coincide with the tableaux in
Example 3.3 in {\color{black}{\cite{FLT2}}}, denoted the same way,
up to a shuffling of the coloured monochromatic numbers. Using Theorem
3.5 in \cite{FLT2}, the corresponding coefficients $ a_{\s} $ are calculated in Example 3.6 in \cite{FLT2}
and are
\begin{equation}\label{as}
  a_{\s_0} = 1, \,  a_{\s_1} = \frac{1}{2}, \, a_{\s_2} = \frac{1}{5}, \,
  a_{\s_3} = -\frac{1}{5}, \, a_{\s_4} = \frac{1}{10}, \, 
  a_{\s_5} = -\frac{1}{10}. 
\end{equation}
(Note that we have here corrected the
values of 
$ a_{\s_2} $ and $  a_{\s_3} $ that are indicated wrongly in \cite{FLT2}). One sees that upon
specializing $ q= 1 $, the coefficients in \eqref{sixexample} become the coefficients in 
\eqref{as}. 

\medskip
We finally remark that in \cite{FLT2} the approach to Theorem {\ref{f_n}} is the converse of
ours, in the sense that it is based on the
characterization of the tableaux in $ \mathcal{CS}_{\lambda} $, whereas in our approach the 
characterization of the tableux in $ \mathcal{CS}_{\lambda} $ is a consequence of Theorem {\ref{f_n}}.
The approach in \cite{FLT2} allows the authors to determine the greatest common divisor of the coefficients
of the expansion of $ f_n$
}

\section{Restricted Specht modules and $ f_\T$ for $\color{black}{ \T}$ a general {\color{black} standard}
  tableau}\label{restricted}
In this section we give an application of the results and methods of the previous sections
to the modular representation theory of the Hecke algebra, that is
the representation theory of $  {\cal H}_n^{  {\Bbbk}}(q)   $, where
$  {\Bbbk} $ is a field made into an $ \cal A $-algebra via $ q \mapsto \xi \in  {\Bbbk}^{\times}$.
We also study the problem 
of expanding $ f_\T $ in terms of standard elements, for $ \T $ a general
{\color{black} standard}
$ \lambda$-tableau
which is not necessarily
a generalized James-Murphy tableau. In this general case we are unfortunately not able
to produce an expansion in terms of standard elements, but still we get some interesting results.

\medskip
There is a natural embedding $  {\cal H}^{\cal A}_{n-1}(q)  \subset  {\cal H}^{\cal A}_n(q) $ which gives rise 
to a restriction functor $ \Res   $ from the category of $ {\cal H}^{\cal A}_n(q)$-modules to the category
of ${\cal H}^{\cal A}_{n-1}(q) $-modules. 

\medskip
In particular, 
for $ \lambda \in \Par_n $ we obtain an $  {\cal H}^{\cal A}_{n-1}(q) $-module $ \Res \, S^{\cal A}_q({\lambda}) $. 
By the {\it branching rule}, it is known that $ \Res \, S^{\cal A}_q({\lambda}) $ 
has a {\it Specht filtration}, that is 
an $  {\cal H}^{\cal A}_{n-1}(q) $-module filtration in which the subquotients are
$  {\cal H}^{\cal A}_{n-1}(q) $
Specht modules.
Let us explain a combinatorial construction of this filtration.

\medskip
Let $ [\alpha_i, \beta_i], i = 0, \ldots, N $ be {\it all} the removable nodes of $ \lambda$,
read from top to bottom as for example in \eqref{fattableauxC}. 
For $ i =0,\ldots, N $ we define $ E_{i} \subseteq S^{\cal A}_q({\lambda})$ via 
\begin{equation}
  E_{i} := \Span_{\cal A} \{ x_\T \, | \, \T \in \std(\lambda)   \mbox{ and } n 
\mbox{ appears in or below the $ \alpha_i ${\color{black}{'}}th row of } \T \,  \}. 
\end{equation}
Then we have that
\begin{equation}\label{isafilttration}
 0 \subset  E_{N} \subset E_{N-1} \subset \cdots \subset 
 E_0 =  \Res \, S^{\cal A}_q({\lambda})
\end{equation}
is a filtration of
$ {\cal H}^{\cal A}_{n-1}(q) $-modules 
and that 
\begin{equation}
  E_{i}/ E_{i+1} \cong S^{\cal A}_q(\lambda^{i}) \mbox{ where } \, \,  \lambda^i:= \lambda \setminus [\alpha_i, \beta_i].
 \end{equation} 
This is a well-known fact, that relies on the Garnir relations, that is
\eqref{Garnir}.
From this filtration we get in particular that $ \Res S_q^{\cal A}({\lambda}) $ contains 
$  S_q^{\cal A}({\lambda^N}) $ as a submodule and has $  S_q^{\cal A}({\lambda^0}) $ as a quotient module.

\medskip
Since the $ E_i $'s and $ S^{\cal A}_q(\lambda^i) $'s are free over $ \cal A $
there is a similar filtration for the specialized Specht module $ \Res  S_q^{ \Bbbk}(\lambda) $, 
which in particular
contains 
$  S_q^{ \Bbbk}({\lambda^N}) $ as a submodule and has $  S_q^{ \Bbbk}({\lambda^0}) $ as a quotient module.

\medskip
Let us now set up some further notation.
Suppose that $ \Bbbk $ is a field which is made into an $ \cal A$-algebra
via $ q \mapsto \xi \in  \Bbbk^{\times} $.
Let $  {\cal A}_{\mathfrak m} $ be 
the localization of $ \cal A $ at the maximal ideal $ {\mathfrak m } := {\ker }({\cal A} \rightarrow \Bbbk)$. 
Then we have that $  {\cal A}_{\mathfrak m}  \subseteq \cal K $. 
Let ${\cal H }^{{\cal A}_{\mathfrak m}}_{n}(q) $ be the Hecke algebra defined over 
$ {\cal A}_{\mathfrak m} $ instead of $ \cal A $.
All constructions for ${\cal H }^{{\cal A}  }_{n}(q) $ can also be carried out for
${\cal H }^{{\cal A}_{\mathfrak m}  }_{n}(q) $ and in particular we have Specht modules
$ S_q^{{\cal A}_{\mathfrak m} }(\lambda)$ for ${\cal H }^{{\cal A}_{\mathfrak m}  }_{n}(q) $. Note
that $ S_q^{{\cal A} }(\lambda) \subseteq 
S_q^{{\cal A}_{\mathfrak m} }(\lambda)   \subseteq S_q^{\cal K}(\lambda) $ and that 
$ S_q^{{\cal A}_{\mathfrak m} }(\lambda) \otimes_{\cal A_{\mathfrak m}} \Bbbk  = 
S_q^{{\cal A} }(\lambda) \otimes_{\cal A} \Bbbk = S_q^{ \Bbbk}({\lambda})  $.

\medskip
As in the previous sections we let $f_n  \in S_q^{\cal K}(\lambda)$
be the seminormal basis element, corresponding to the James-Murphy tableau
$ \T_{n} $, for example as in \eqref{exJMtableau}. 
We have $  {\cal H }^{{\cal A}  }_n(q)
\subseteq  {\cal H }^{\cal K  }_n(q) $ and so we may introduce an 
$   {\cal H }^{{\cal A}  }_{n-1}(q) $-module  $ U_q^{\cal A}({\lambda^0})  $ as follows
\begin{equation} U_q^{\cal A}({\lambda^0}) := 
f_n   {\cal H }^{{\cal A}  }_{n-1}(q)
\subseteq  S_q^{\cal K}(\lambda).
\end{equation}
Assume now further that $ [r_{d_i}]_{\xi} \not= 0 $ for $ i = 1, 2, \ldots , N $.
Then by Theorem {\ref{f_n}} we have that $ f_n \in  S_q^{{\cal A}_{\mathfrak m} }(\lambda)$
and so we can define an $ {\cal H }^{{\cal A}_{\mathfrak m}  }_{n-1}(q) $-submodule
$ U_q^{{\cal A}_{\mathfrak m}}({\lambda^0})  $ 
of
$ \Res S_q^{{\cal A}_{\mathfrak m}}(\lambda) $ via 
\begin{equation} 
  U_q^{{\cal A}_{\mathfrak m}}({\lambda^0})
  := f_n  {\cal H }^{{\cal A}_{\mathfrak m}  }_{n-1}(q)
\subseteq  \Res S_q^{{\cal A}_{\mathfrak m}}(\lambda).
\end{equation}

{\color{black}Note that $ U_q^{{\cal A}_{\mathfrak m}}({\lambda^0}) $ is not defined 
  as a specialization of $ U_q^{{\cal A}}({\lambda^0}) $ since that would not be
a submodule of $ \Res S_q^{{\cal A}_{\mathfrak m}}(\lambda) $}. 

\medskip
We now have the following Theorem.

\begin{Thm}{\label{modular_Specht}} 
\begin{enumerate}[label=\alph*)]  \itemsep 0.2cm
\item There is an $   {\cal H }^{{\cal A}  }_{n-1}(q) $-isomorphism 
$S_q^{{\cal A}  }{(\lambda^0)} \rightarrow U_q^{\cal A}(\lambda^0)    $  
  given by $x_{\lambda^0} \mapsto f_n   $ and similarly, when
    $ [r_{d_i}]_{\xi} \not= 0 $ for $ i = 1, 2, \ldots , N $, there is an 
  $   {\cal H }^{{\cal A}_{\mathfrak m}  }_{n-1}(q) $-isomorphism
  $S_q^{{\cal A}_{\mathfrak m}  }{(\lambda^0)} \rightarrow U_q^{{\cal A}_{\mathfrak m}}(\lambda^0)    $
    given by $x_{\lambda^0} \mapsto f_n   $.
\item   Suppose that 
  $ [r_{d_i}]_{\xi} \not= 0 $ for $ i = 1, 2, \ldots , N $.
  Then $  S_q^{ \Bbbk}({\lambda^0}) $ splits off from 
  $ \Res S_q^{\Bbbk  }({\lambda})  $ with splitting homomorphism given by specializing
  the composition
  $S_q^{{\cal A}_{\mathfrak m}  }{(\lambda^0)}  \rightarrow U_q^{{\cal A}_{\mathfrak m}}(\lambda^0)  \subseteq
   \Res S_q^{{\cal A}_{\mathfrak m} }({\lambda})   $.
  \end{enumerate}
\end{Thm}
\begin{dem}
To show $a) $ 
we first verify that $ f_n \mapsto x_{\lambda^0} $ defines an
$   {\cal H }^{{\cal A}  }_{n-1}(q) $-homomorphism  
$  U_q^{\cal A}(\lambda^0)  \rightarrow S_q^{{\cal A}  }{(\lambda^0)} $. 
This is not completely obvious, since
$ f_n $ and $  x_{\lambda^0} $ may apriori have different annihilators in
$   {\cal H }^{{\cal A}  }_{n-1}(q) $. We resolve this problem as follows. 
By definition $ U_q^{\cal A}(\lambda^0) $ is an $   {\cal H }^{{\cal A}  }_{n-1}(q) $-submodule
of the $   {\cal H }^{{\cal K}  }_{n-1}(q) $-module 
$ U_q^{\cal K}(\lambda^0)  := f_n {\cal H }^{\cal K}_{ n-1}(q) $. 
Similarly, the Specht module $ S_q^{{\cal A}  }{(\lambda^0)} = x_{\lambda^0}  {\cal H }^{{\cal A}  }_{n-1}(q) $
is an $   {\cal H }^{{\cal A}  }_{n-1}(q) $-submodule 
of $ S_q^{{\cal K}  }{(\lambda^0)} =  x_{\lambda^0}  {\cal H }^{{\cal K}  }_{n-1}(q)  $. 
Define 
\begin{equation}
\std^0(\lambda) := \{ \T \in \std(\lambda) \mid \T[\alpha_0, \beta_0] = n\}. 
\end{equation}
Using YSF,  that is Theorem {\ref{Seminormal_form}}, we have that $ U_q^{\cal K}(\lambda^0) $ is
generated by 
$
\{f_{\T} \, | \, \T \in \std^0(\lambda)  \} 
$
and since $ \{ f_\T  \,| \, \T \in \std^0(\lambda) \}   \subseteq
 \{ f_\T \, | \, \T \in \std(\lambda) \}   \subseteq S_q^{\cal K}(\lambda) $ 
we also have that $
\{f_{\T} \, | \, \T \in \std^0(\lambda)  \} 
$
is $ \cal K $-linearly independent. Hence it is a $ \cal K $-basis for $ U_q^{\cal K}(\lambda^0)  $. On the
other hand, 
$ \{f_{\s} \, |\, \s \in \std(\lambda^0)  \}  $ is a $ \cal K $-basis for $ S^{\cal K}_q(\lambda^0)  $,
and so we obtain a $ \cal K $-linear bijection 
$ \varphi: 
U_q^{\cal K}(\lambda^0)  \rightarrow S_q^{{\cal K}  }{(\lambda^0)}$, via 
\begin{equation}
 \varphi( f_{\T} ) = f_{\s} \mbox{ where } \T \in  \std^0(\lambda) \mbox{ and }   \s=  \T |_{ 1,2, \ldots, n-1 } 
\end{equation}
Using YSF on $f_{\T}  $ as well as on $ f_{\s} $, we conclude that $ \varphi $ is in fact 
an ${\cal H }^{{\cal K}  }_{n-1}(q) $-isomorphism, since the action of $ T_i $ on both cases is
'the same', given only by radial lengths.
The restriction of $ \varphi $ to
$ U_q^{\cal A}(\lambda^0) $ is the inverse of
the ${\cal H }^{{\cal A}  }_{n-1}(q) $-isomorphism that is postulated in $ a)$.
The second part of $a) $, involving the ground ring $ {\cal A}_{\mathfrak m} $, is proved
the same way.

\medskip
To show $b)$, letting $ \psi: S_q^{{\cal A_{\mathfrak m}}  }{(\lambda^0)} \rightarrow
U_q^{{\cal A}_{\mathfrak m}}(\lambda^0) $ be the isomorphism from $a)$, we have that
$ \psi(x_{\lambda^0}) =    f_n   $.  Since $ U_q^{{\cal A}_{\mathfrak m}}(\lambda^0) \subseteq 
 \Res S_q^{{\cal A}_{\mathfrak m}  }({\lambda})  $ we obtain by specializing an 
$   {\cal H }^{\Bbbk }_{n-1}(q) $-homomorphism $ (\iota \circ \psi) \otimes 1 :
 S_q^{\Bbbk  }{(\lambda^0)} \rightarrow \Res S_q^{\Bbbk  }({\lambda})  $,
where $ \iota: U_q^{{\cal A}_{\mathfrak m}}(\lambda^0) \rightarrow 
 \Res S_q^{{\cal A}_{\mathfrak m}  }({\lambda}) $ is the inclusion homomorphism. Note that 
 apriori $ (\iota \circ \psi) \otimes 1 $ may not be injective, since the functor $ ? \otimes \Bbbk$
 is not left exact. But letting 
 $ \pi \otimes 1: \Res S_q^{\Bbbk  }({\lambda}) \rightarrow   S_q^{ \Bbbk}({\lambda^0}) $
be the quotient map from the specialization of the  
filtration \eqref{isafilttration}, 
we get that $ ((\pi \otimes 1)  \circ ( \iota \circ \psi) \otimes 1)( x_{\lambda^0} \otimes 1)  =  x_{\lambda^0}\otimes 1  $,
and hence $ ((\pi \otimes 1)  \circ ( \iota \circ \psi) \otimes 1) $ is the
identity map on $ S_q^{ \Bbbk}({\lambda^0}) $.
This proves $b)$.
\end{dem}

\medskip
Our next Theorem 
establishes a converse of $ b) $ of the previous Theorem, but over $ {\cal A}_{\mathfrak m}$
instead of $ \Bbbk$.
Its formulation was 
influenced by Proposition 3.11 of
\cite{FLT1}. 
Note however that the authors of \cite{FLT1} work in the setting of the
induced Specht module,
rather than the restricted Specht module as we do. 
Note also {\color{black}{that}}
%{\color{green}{\sout{they}}}
their argument, unlike ours, 
does not rely on the idempotents
$ E_{\T}$. We believe that our argument is more {\color{black}{conceptual}} than the one of \cite{FLT1}.
\begin{Thm}\label{splittingcriteria}
  Let $ \pi $ be the
$   {\cal H }^{{\cal A}_{\mathfrak m}  }_{n-1}(q) $-quotient map $ \pi:  \Res S_q^{{\cal A}_{\mathfrak m}  }({\lambda})
\rightarrow S_q^{{\cal A}_{\mathfrak m}} ( \lambda^{0}) $. Then $ \pi $ admits a splitting if and only if
$ [r_{d_i}]_{\xi} \not= 0 $ for $ i = 1, 2, \ldots , N $. 
\end{Thm}
\begin{dem}
  Suppose that $ [r_{d_i}]_{\xi} \not= 0 $ for $ i = 1, 2, \ldots , N $. Then we actually showed in 
  $ b) $ of Theorem {\eqref{modular_Specht}} 
  that $ \pi $ has a splitting as an $   {\cal H }^{{\cal A}_{\mathfrak m}  }_{n-1}(q) $-homomorphism
  and so we only need to prove the 'only if' part of the Theorem.

\medskip  
Suppose therefore that $ \psi: S_q^{{\cal A}_{\mathfrak m}} ( \lambda^{0}) \rightarrow
\Res S_q^{{\cal A}_{\mathfrak m}  }({\lambda}) $ is a splitting homomorphism for $ \pi$ and set
$ f_n^{\prime} := \psi(x_{\lambda^{0}})$. 
Then $ f_n^{\prime} := \psi(x_{\lambda^0} ) $ has an expansion 
\begin{equation}\label{fnhasanexpanson}
f_n^{\prime} = x_{\T_n} + \sum_{\s \in \std(\lambda)} a_{\s} x_\s, \, a_{\s} \in {\cal A}_{\mathfrak m}. 
\end{equation}  
By the splitting property $  \pi \circ \psi = Id $, we have that 
for all $ \s $ occurring in the sum, $ n $ appears in $ \s $ {\color{black}in} a node strictly below
$ [ \alpha_0, \beta_0]$. 

\medskip
Extending scalars from $ {\cal A}_{\mathfrak m} $ to $ \cal K $ we now consider $ f_n^{\prime}  $
as an element of $ \Res S_q^{\cal K  }({\lambda}) $ and similarly we consider $ \pi $ and $ \psi $ 
as $  {\cal H }^{\cal K  }_{n-1}(q) $-homomorphisms. Suppose now that
$ \T \in \std(n-1) $ and let 
$ E_{\T} \in   {\cal H }^{\cal K  }_{n-1}(q) $
be the idempotent corresponding to $ \T$, as introduced in  
{\eqref{E-expansion}}. We then have that 
\begin{equation}\label{combining1}
  f_n^{\prime}   E_{\T} =  \psi ( x_{\lambda^0})  E_{\T} =
  \psi ( x_{\lambda^0}E_{\T})   = \left\{ \begin{array}{ll} \psi ( x_{\lambda^0}) =   f_n^{\prime}  & \mbox{ if }
      \T = \T^{\lambda^0} \\ 0 & \mbox{ otherwise.} \end{array} \right.
\end{equation}  
We are interested in the idempotent $ E_{\T_n} \in {\cal H }^{\cal K  }_{n}(q) $
since we must show that
$x_{\T_n}  E_{\T_n} =   f_n^{\prime}$.

\medskip
\noindent
Note that for any $ \T \in \std(n) $ 
we have the following formula which can be read off 
from {\eqref{E-expansion}} 
\begin{equation}\label{combining2}
 E_{\T} = E_{\T|_{ 1,2, \ldots, n-1 }} E_{\T|_{ n }} \mbox{ where }  E_{\T|_{ n }}:= 
 \prod_{c \in \Cont_n \setminus c_\T(n) }
\, \frac{L_n -c }{ c_\T(n) -c}.
\end{equation}  
For all $ \s $ occurring in \eqref{fnhasanexpanson}
we have that
\begin{equation}\label{combining0} x_\s  E_{\T_n}=0. \end{equation}
Indeed, we have the general triangularity property
\begin{equation}\label{combining3}
x_{\U} E_{\V}  \neq 0 \Longrightarrow \U \unlhd \V 
\end{equation}
which follows from inverting the expansion \eqref{ftriangular} and using
the orthogonality of the $ E_{\V}$'s, and so if $ x_\s  E_{\T_n} $ were nonzero, we would have
$ \s \unlhd \T_n$. But removing the $n$-node from both sides we have
$  shape(\s|_{ 1,2, \ldots, n-1 })  \rhd shape(\T_n|_{ 1,2, \ldots, n-1 }) = \lambda $, in contradiction with
$ \s \unlhd \T_n$.

\medskip
Similarly to \eqref{combining3} we have that
\begin{equation}\label{combining4}
x_{\U} E_{\V}  \neq 0 \Longrightarrow shape(\U) = shape(\V).
\end{equation}

Combining \eqref{combining0}, \eqref{combining2}, \eqref{combining3} and \eqref{combining4}
and using that $ \sum_{\T \in \std(n)}  E_{\T} =1 $ we now get
\begin{equation}
  f_n^{\prime} = f_n^{\prime} \sum_{ t\in \std(n)} E_{ \T } =
  f_n^{\prime} \sum_{ t\in \std(\lambda)} E_{ \T } =
    f_n^{\prime} E_{\T_n} =   
  x_{\T_n} E_{\T_n} = f_n
\end{equation}
which proves the Theorem.
\end{dem}

\medskip
\noindent
{\bf Remark}
For $M, N$ general $ {\cal H }^{\cal A  }_{n}(q) $-modules it is not true
that any $  {\cal H }^{\Bbbk  }_{n}(q) $-homomorphism $ \varphi $ between the specialized modules
$ M \otimes \Bbbk $ and $ N \otimes \Bbbk $ can be lifted to an $  {\cal H }^{\cal A  }_{n}(q) $-homomorphism
between $ M $ and $ N $.
On the other hand, 
if $M $ and $N$ are $ {\cal H }^{{\cal A}_{\mathfrak m}  }_{n}(q) $-modules, then it appears plausible that 
$ \varphi: M \otimes \Bbbk \rightarrow N \otimes \Bbbk  $ can be lifted. This issue 
is discussed in \cite{FLT1} for the splitting of the top Specht quotient from the induced Specht module  
in
{\color{black}{the}}
symmetric group case, where in particular it is pointed out 
that there are no known counterexamples to the lifting property 
in this case.
If the lifting property hold{\color{black}s} for 
$ {\cal H }^{{\cal A}_{\mathfrak m}  }_{n}(q) $-modules, then
we could improve the criterion of Theorem \ref{splittingcriteria} to a criterion for 
$ {\cal H }^{\Bbbk  }_{n}(q) $-splittings.

\medskip
We now return to our main object of study, namely that of expanding $ f_\T$ in terms
of standard elements.

\medskip
We can generalize $ a) $ of Theorem {\ref{modular_Specht}} as follows. 
Let $ \T $ be a {\color{black} standard} $\lambda$-tableau.
Then $ \T $ can be viewed as a chain of
partitions $ \{ \lambda^{\le j} \}_{j=1,\ldots, n} $ where $ \lambda^{\le j} :=
shape(\T|_{ 1,2, \ldots, j  }) $. 
For $ 1 \le j \le n  $ we define the $ \lambda $-tableau
$ \T^{\leq j }_{ | j+1, \ldots, n} $ via
\begin{equation}{\label{ex103}}
  ( \T^{\leq j }_{ | j+1, \ldots, n})^{-1}(k) :=
  \left\{ \begin{array}{ll} (\T^{\lambda^{ \leq j}})^{-1}(k) & \mbox{ if } k \le j \\
      \T^{-1}(k) & \mbox{ if } k > j. \end{array} \right.
\end{equation}

For simplicity we also write $ \T^{\leq j}= \T^{\leq j }_{ | {\color{black}{j+1, \ldots, n}}} $.
{\color{black}Below is an example of a tableau $ \T $ with corresponding tableau
  $ \T^{\leq 10} $. The numbers $ k \le 10 $, corresponding to the first case in
  \eqref{ex103}, have been coloured blue. 
\begin{equation}
  \T = \raisebox{-.5\height}{\includegraphics[scale=0.7]{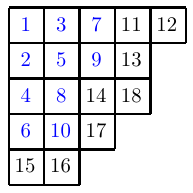}}, \, \, \, \, \, \, 
   \T^{\leq 10} = \raisebox{-.5\height}{\includegraphics[scale=0.7]{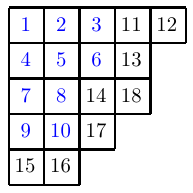}} \,  
\end{equation}}

We then introduce the $ {\cal H }_j^{ \cal A}(q) $-submodule $  U_q^{\cal A}({ \T^{\leq j} }) $ 
of $ S_q^{\cal K}({ \lambda}) $ via 
$ U_q^{\cal A}({ \T^{\leq j} })
:= f_{ \T^{ \leq j} }  {\cal H }_j^{ \cal A}(q) $. 
This is a generalization of $U_q^{\cal A}({\lambda^0}) $ from
Theorem {\ref{modular_Specht}} which is recovered by setting 
$ j := n -1$. 
We now have the following Theorem,
generalizing $a)$ of Theorem \ref{modular_Specht}.

\begin{Thm}{\label{modular_Specht_1}} 
  There is an isomorphism of $ {\cal H }^{ \cal A}_j(q)$-modules 
  \begin{equation}  \varphi: 
 U_q^{\cal A}({ \T^{\leq j} })
    \rightarrow 
S^{ { \cal A} }_q({ \lambda^{ \le j}}), 
\, \, \, \, \, \, \, \, \, \, \,  
f_{ t^{ \leq j} } 
\mapsto  x_{ \lambda^{ \le  j } }. 
\end{equation}
\end{Thm}
\begin{dem}
The proof is essentially the same as the proof of $a)$ of the previous Theorem. 
Let us briefly indicate 
the necessary modifications.
By definition $  U_q^{\cal A}({ \T^{\leq j} }) $ is an
$ {\cal H }^{ \cal A}_j(q)$-submodule of the $  {\cal H }_j^{ \cal K}(q)$-module
$ U_q^{\cal K}(  \T^{\leq j}    )  := f_{ \T^{ \leq j} } {\cal H }^{\cal K}_{ n-1}(q) \subseteq S_q^{\cal K}(\lambda) $.
Using YSF, we see that $ U_q^{\cal K}(  \T^{\leq j}    )  $ has $ \cal K $-basis
\begin{equation}
\{  f_\s \, | \, \T \in \std(\lambda) \mbox{ and } \s^{-1}(k) = \T^{-1}(k)  \mbox{ for } k =j+1, \ldots, n \}.
\end{equation}
But then, using the $ \cal K$-basis 
$ \{  f_\U \, | \, \U \in \std(\lambda^{ \le  j }) \}$
for $ S_q^{\cal K}( \lambda^{ \le  j }) $, we get via YSF that $ f_\s \mapsto f_{\s_{  | 1, \ldots, j}} $ induces
an $ {\cal H }^{\cal K}_{ j}(q) $-isomorphism $   U_q^{\cal K}(  \T^{\leq j}    )  \rightarrow 
S_q^{\cal A}( \lambda^{ \le  j }) $. As the required $ \varphi $ we can then use
the restriction of this isomorphism to
$ U^{\cal A}_q(  \T^{\leq j}    ) $. 
\end{dem}

\medskip
With this result established, we now 
finally consider the problem of determining the expansion in standard elements of
$ f_\T $ where 
this time $ \T $ is a completely general {\color{black} standard} $\lambda$-tableau.

\medskip
Given $ \T$, 
we define an element $ {\cal P}_{\color{black}{\T}} \in   {\cal H}_n    $ via the following formula
\begin{equation}\label{generalexpansion}
  {\cal P}_{\color{black}{\T}} :={\cal  P}_{n}\,  {\cal P}_{n-1}  \cdots    {\cal P}_{1	}
\end{equation}  
where 
$ {\cal P}_j  \in {\cal H}_j \subseteq {\cal H}_n $ is
the element $ {\cal P}^N$ given by the recursion \eqref{mainrecursion}, but with respect to
the $ \lambda^{\le j }$-tableau
$ \T^{\leq j-1 }_{ | j} $.

\medskip
We now have the following generalization of 
Theorem {\ref{f_n}}.
\begin{Thm}{\label{final}} 
In the above setup we have $ f_\T =  x_{\lambda} \, {\cal P}_\T $.
\end{Thm}
\begin{dem}
By Theorem \ref{f_n}
we have that $ f_{\T_n} = f_{\T_{| n}^{ \leq n-1 }}= x_{\lambda} \, {\cal P}_n $
and also 
$ f_{\T_{| n-1}^{\leq n-2 }} =  x_{ \lambda^{\leq (n-1)} } {\cal P}_{n-1} $. In view of Theorem
{\eqref{modular_Specht_1}} we deduce from this that 
\begin{equation}
f_{\T_{| n-1, n}^{\leq n-2 }} = f_{\T_{| n}^{\leq n-1 }} {\cal P}_{n-1} = 
x_{\lambda} \, {\cal P}_n  {\cal P}_{n-1}. 
\end{equation}  
This argument is now repeated until arriving at 
the formula claimed in the Theorem.
\end{dem}

\medskip
Let us illustrate the Theorem on the partition $ \lambda = (3, 1^2 ) $ of $5 $ and 
the $\lambda$-tableau
\begin{equation}\label{illustrativeexample}
\T:=   \raisebox{-.5\height}{\includegraphics[scale=0.7]{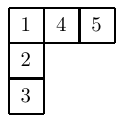}}. 
\end{equation}

\medskip
We have by the Theorem that $ f_\T = x_{\lambda} {\cal P}_5 {\cal P}_4 {\cal P}_3 {\cal P}_2 {\cal P}_1 $.
Let us work out the ${\cal P}_i$'s. 
Since $ \T_{| 1,2,3} $ is the largest, in fact the only,  standard $ \lambda^{\leq 3} $-tableau
we find immediately that $ {\cal P}_3 = {\cal P}_2 = {\cal P}_1 = 1 $.
Let us then consider $ {\cal P}_5$. 
By Theorem \ref{fat_hook_lemma} we have for $ f_{\T^{\leq 4}_{| 5}} $ 
the following expansion in terms of standard $ x_\T $'s. 
\begin{equation}\label{app}
  f_{t^{\leq 4}_{| 5}} \, \, = \, \,
\raisebox{-.5\height}{\includegraphics[scale=0.7]{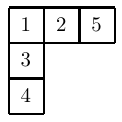}} 
+ \, \, \frac{1}{\, [4]_q} \left( \,
\raisebox{-.5\height}{\includegraphics[scale=0.7]{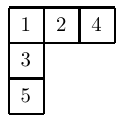}} 
-q \raisebox{-.5\height}{\includegraphics[scale=0.7]{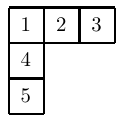}} \right). 
\end{equation}

Thus we have $ {\cal P}_5 = T_{3,5} + \frac{1}{\, [4]_q} (T_3  -q ) $.
We next work out $ {\cal P}_4 $. Using Theorem \ref{fat_hook_lemma} once 
more we have 
for $ f_{\T^{\leq 3}_{| 4}} $ 	the following 
expansion

\begin{equation}
  f_{\T^{\leq 3}_{| 4}} \, \, = \, \,
\raisebox{-.5\height}{\includegraphics[scale=0.7]{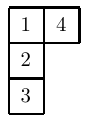}} 
+ \, \, \frac{1}{\, [3]_q} \left( \,
\raisebox{-.5\height}{\includegraphics[scale=0.7]{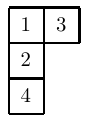}} 
-q \raisebox{-.5\height}{\includegraphics[scale=0.7]{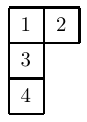}} \right)
\end{equation}
and so we 
get $ {\cal P}_4 = T_{2,4}  + \frac{1}{\, [3]_q} ( T_2-q) $. 
Combining, we get an expression for $ f_\T = x_\lambda {\cal P}_5 {\cal P}_4 {\cal P}_3 {\cal P}_3 {\cal P}_3
= x_\lambda {\cal P}_5 {\cal P}_4 $
in terms of
9 $ x_\T $'s, which however unfortunately are not all standard. After straightening we get 
\begin{equation} f_\T = 
\raisebox{-.5\height}{\includegraphics[scale=0.7]{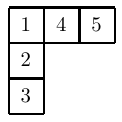}} 
+ \, \, \frac{1}{\, [3]_q} \left( \,
\raisebox{-.5\height}{\includegraphics[scale=0.7]{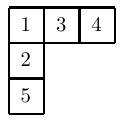}} 
-q \raisebox{-.5\height}{\includegraphics[scale=0.7]{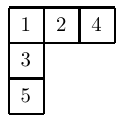}}
-q \raisebox{-.5\height}{\includegraphics[scale=0.7]{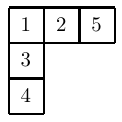}}
- \raisebox{-.5\height}{\includegraphics[scale=0.7]{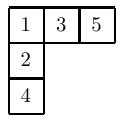}}
\right).
\end{equation}

For comparison, repeated use of 
Young's seminormal form twice, that is the algorithm explained in section \ref{analgorithm} on 
$ f_\T $ would have given $ 4 \cdot 3 = 12$ $ x_\T$'s 
instead of $9$ (that after straightening would have reduced to the above expression, of course).
In general, as actually already follows from Theorem {\ref{f_n}}, this algorithm will
in {\color{black}{general}} produce more 
than just one denominator, i.e. the above example with the only denominator $[3]_q $ is special in this
respect.

\medskip
\noindent
    {\bf Remark} In general, as we just saw on the example \eqref{illustrativeexample},
    the expansion of $ f_\T$ using Theorem {\ref{final}} does not always produce standard
    elements $x_\T$. It is an interesting open problem to find an {\color{black}
      efficient}} algorithm that does give 
    such an expansion. 

\medskip    
    {\color{black}In spite of this remark, we can still use Theorem \ref{final} to deduce
      the following consequence for the coefficient of $ x_\s$ of $ f_\T $, which is valid for 
      general standard tableaux $ \s $ and $ \T$. We are grateful to the referee for pointing this out to
      us.
      \begin{Cor}\label{cor}
Let $ \{f_\T \} $ and $ \{ x_\T\}  $ be the seminormal and standard bases for $ {\cal H }^{{\C(q)}}_{ n}(q)   $ 
and let $ f_\T = \sum_{ \s \in \std(\lambda) } c_{\s \T } x_\s $ be the
expansion of $ f_\T $, with $ c_{\s \T } \in   {\mathbb C}(q)$. Then for
any $ \s, \T \in \std(\lambda) $, the poles of 
$ c_{\s \T } $ are roots of unity in $ \C$.  
\end{Cor}       
\begin{dem}
  The poles of $  \dfrac{1}{[k]_q} = \dfrac{q-1}{q^k -1} $ are roots of unity
  and so via Theorem {\ref{final}}, together with the definitions in
\eqref{mainrecursion} and \eqref{generalexpansion}, 
we get that the coefficient of $ x_\s $
in the expansion of $ f_\T $, where $ \s \in \tab(\lambda) $, 
has poles that are roots of unity.
But the expansion of $ x_\s $ in terms of standard elements, using the Garnir relations, can be
carried out over $ \cal A $ and will therefore not introduce new poles. This proves the Corollary.
  \end{dem} 

}  
\medskip
Let us do a rudimentary complexity analysis of the two algorithms for calculating 
$ f_\T$, that is 'repeated use of 
Young's seminormal form' versus the algorithm given by Theorem {\ref{final}}. 
Suppose $ \lambda = (\lambda_1 , \lambda_2^{k_2 } ) $ is a fat hook partition with {\color{black}{the}}
first block of rows of width one, that is $ \lambda_1 > \lambda_2$.
We have $ n = \lambda_1  + \lambda_2 k_2 $ and $ a = \lambda_1  $, that is 
$ \sigma_{a,n} $
%{\color{green}{\sout{has}}}
is of Coxeter length $ \lambda_2 k_2   $. Thus, the repeated use of Young's 
seminormal
{\color{black}{form}}
to calculate $ f_n $ produces a linear combination $ 2^{ \lambda_2 k_2} $
(standard and nonstandard) elements
$x_\T$'s, whereas the algorithm contained in Theorem \ref{fat_hook_lemma} produces
%for example in the 
%formulation given in Corollary 
%{\ref{hook_lemma_cor}},   
\begin{equation} (\lambda_2 -1 )^{k_2 }  + \ldots + (\lambda_2 -1)^2 + (\lambda_2 -1) = 
  \frac{(\lambda_2 -1 )^{k_2 +1 } - \lambda_2 -1}{\lambda_2  -2 }
\end{equation} 
such elements $ x_\T$. Thus, with respect to $ \lambda_2  $ we see that  
Theorem \ref{fat_hook_lemma} has polynomial complexity whereas 
repeated use of Young's seminormal form 
has exponential complexity. Thus the algorithm contained in Theorem \ref{fat_hook_lemma}
is much more efficient. 
This relationship carries over to the general 
algorithm of Theorem {\ref{final}}. 
We have implemented the algorithms using the GAP system.

\medskip
\noindent
{\bf Remark}
As already mentioned there is no known algorithm for
expanding a general $ f_\T $ in terms of standard elements. On the other hand, our calculations
for small (but nontrivial) $ \T $, using the Gram-Schmidt algorithm explained in \eqref{downwards},
indicate that the coefficients of the expansion are 'nice' expressions involving radial lengths.
Unfortunately, we are at this point unable to state the exact meaning of 'nice'.

steen@inst-mat.utalca.cl, steenrh@gmail.com, Universidad de Talca, Chile.


\bigskip
\begin{thebibliography}{X}
\bibitem[Ar]{Ar} S. Ariki, \textit{On the decomposition numbers of the Hecke algebra
$ G(m,1,n) $}, J. Math. Kyoto Univ. {\bf 36}, (1996), 789-808. 

{\color{black}{
\bibitem[AH]{AH}
S. Armon, T. Halverson, \textit{
Transition matrices between Young’s natural and seminormal representations}, 
Electron. J. Comb. {\bf 28}, (2021), No. 3, Research Paper P3.15, 34 p.
}}

\bibitem[C]{C} I. Cherednik, \textit{Double affine Hecke algebras and Macdonald's conjectures}, Annals of
Maths., {\bf 141}, (1995), 191-216. 
\bibitem[DJ]{DJ} R. Dipper, G. James,  \textit{Blocks and idempotents of Hecke algebras of general linear 
groups},  Proc. London Math. Soc. (3) {\bf 54}, (1987), no. 1, 57--82. 

\bibitem[GL]{GL} J. J. Graham, G. I. Lehrer,
\textit{Cellular algebras}, Inventiones Mathematicae {\bf 123}, (1996), 1-34.


\bibitem[H]{H} M. Haiman, \textit{Hilbert schemes, polygraphs and the Macdonald positivity conjecture}, 
J. Amer. Math. Soc. {\bf 14}, (2001), 
no. 4, 941-1006. 

\bibitem[J]{J} G. D. James, \textit{The representation theory of the 
symmetric groups}, Lecture notes in mathematics {\bf 682}, Springer Verlag
(1978).
\bibitem[JM]{JM} G. D. James, G. E. Murphy, \textit{The Determinant of the
Gram Matrix for a Specht Module}, Journal of Algebra {\bf 59}, (1979), 222-235. 
\bibitem[LLT]{LLT} A. Lascoux, B. Leclerc, J.-Y. Thibon,
{\textit Hecke algebras at roots of unity and crystal bases of quantum affine 
algebras}, Commun. Math. Phys. {\bf 181}, (1996), 205-263.
\bibitem[FLT1]{FLT1}
M. Fang, K. J. Lim,  K. M. Tan, \textit{Jantzen filtration of Weyl modules, product of Young
    symmetrizers and denominators of Young’s seminomal basis}, 
Represent. Theory {\bf 24}, (2020), 551-579.
\bibitem[FLT2]{FLT2} M. Fang, K. J. Lim,  K. M. Tan, \textit{Young’s seminormal basis vectors and their denominators}, {\color{black}{Journal of Combinatorial Theory, Series A, 
{\bf 184}, (2021), 105494}}. 









  
\bibitem[M1]{M1} G. E. Murphy, \textit{On the Representation Theory of the 
Symmetric
Groups and Associated Hecke Algebras}, Journal of Algebra {\bf 152}, (1992),
492-513.
\bibitem[M2]{M2} G. E. Murphy, \textit{The Representations of Hecke Algebras of type
$ A_n $}, Journal of Algebra {\bf 173}, (1995), 97-121. 
\bibitem[M3]{M3} G. E. Murphy, \textit{A new construction of Young's seminormal representation of the Symmetric groups}, Journal of Algebra {\bf 69}, (1981), 287-297.



  %\bibitem[M5].G. E. Murphy, Composition factors of Specht modules for
%Hecke algebras of type $ A_n$, 
%Journal of Algebra {\bf 306} (2006), 268-289.
\bibitem[Ma]{Ma} A. Mathas, \textit{Iwahori-Hecke Algebras and Schur Algebras of the Symmetric Group}, Univ. Lecture
Series {\bf 15}, (1999), Amer. Math. Soc.
\bibitem[Ram]{Ram} A. Ram, \textit{Skew representations are irreducible}, 
in  Combinatorial and Geometric representation theory, 
S.-J. Kang and K.-H. Lee eds., Contemp. Math. {\bf 325}, Amer. Math. Soc. (2003), 161-189. 
\bibitem[Rai]{Rai} C. Raicu, \textit{Products of Young symmetrizers and ideals in the generic tensor algebra},
  J. Algebraic Combinatorics {\bf 39}, (2014), 247-270.
\bibitem[RH]{RH} S. Ryom-Hansen, \textit{Grading the
  translation functors in type $A$}, J. Algebra {\bf 274}, (2004), no. 1, 138--163.
\bibitem[S{\etalchar{+}}97]{Sch97}
            Martin Sch{\accent127 o}nert et~al.
  \newblock {{GAP} --
            {Groups}, {Algorithms}, and {Programming} --
            version 3 release 4 patchlevel 4"}.
  \newblock Lehrstuhl D f{\accent127 u}r Mathematik,
            Rheinisch Westf{\accent127 a}lische
            Technische Hoch\-schule, Aachen, Germany, 1997.

\end{thebibliography}
\end{document}